\documentclass[%
  onecolumn
, a4paper
]{mpi2015-cscpreprint}

\usepackage[utf8]{inputenc}
\usepackage{amsmath} 
\usepackage{amsfonts}
\usepackage{amssymb}
\usepackage{amsmath}
\usepackage{pgfplots}
\pgfplotsset{compat=newest}
\usepackage{mathrsfs} 
\usepackage[ruled,linesnumbered]{algorithm2e} 
\SetEndCharOfAlgoLine{}
\usepackage{pgfplotstable}
\usepackage{xifthen}

\usepackage{tikz}
\usetikzlibrary{calc, positioning}
\usetikzlibrary{shapes}
\usetikzlibrary{shapes.misc} 
\usepackage{pgfplots}
\pgfplotsset{compat = 1.13,
  unbounded coords = jump,
}
\usepgfplotslibrary{external}
\tikzset{external/system call = {%
    lualatex \tikzexternalcheckshellescape
    -halt-on-error
    -interaction=batchmode
    -jobname "\image" "\texsource"}}
\tikzexternaldisable{}

\usepackage{filemod}

\usepackage[amsmath,amsthm,thmmarks]{ntheorem}
\usepackage{thmtools}

\newtheorem{assumption}{Assumption}
\newtheorem{theorem}{Theorem}
\crefname{assumption}{Assumption}{Assumptions}
\newcommand{\Gammacool}{\ensuremath{\Gamma_{\text{cool}}}}
\newcommand{\GammaN}{\ensuremath{\Gamma_\text{N}}}
\newcommand{\Gammaint}{\ensuremath{\Gamma_{\text{int}}}}

\newcommand{\Gammain}{\ensuremath{\Gamma_{\control}}}
\newcommand{\domain}{\ensuremath{\Omega}}
\newcommand{\domainl}{\ensuremath{\domain_l}}
\newcommand{\domains}{\ensuremath{\domain_s}}
\newcommand{\statexstar}{\ensuremath{\state^*}}
\newcommand{\Gammaintdiff}{\ensuremath{\Gamma_{\text{int},\Delta}(\ttime,\statexstar)}}

\newcommand{\GammaintRef}{\ensuremath{{\Gamma}_{\text{int,ref}}(\ttime, \statexstar)}}
\newcommand{\GammaintRefx}{\ensuremath{{\Gamma}_{\text{int,ref}}(\ttime, \xx)}}
\newcommand{\GammaintReft}{\ensuremath{{\Gamma}_{\text{int,ref}}(\ttime)}}
\newcommand{\Gammainx}[1][]{%
\ifthenelse{\isempty{#1}}{\Gammain}{\ensuremath{\Gamma_{\control,#1}}}}%
\newcommand{\Gammacoolx}[1][]{%
\ifthenelse{\isempty{#1}}{\Gammacool}{\ensuremath{\Gamma_{\text{cool},#1}}}}%
\newcommand{\Gammaoutx}[1][]{%
\ifthenelse{\isempty{#1}}{\ensuremath{\Gamma_{\cC}}}{\ensuremath{\Gamma_{\cC,#1}}}}%
\newcommand{\Tcool}{\ensuremath{\temp_\text{cool}}}
\newcommand{\Tmelt}{\ensuremath{\temp_\text{melt}}}
\newcommand{\TmeltDiff}{\ensuremath{\temp_{\Delta,\text{melt}}}}
\newcommand{\TmeltRef}{\ensuremath{\tilde \temp_\text{melt}}}
\newcommand{\Tzero}{\ensuremath{\temp_0}}
\newcommand{\kk}{\ensuremath{k}}
\newcommand{\ks}{\ensuremath{\kk_s}}
\newcommand{\kl}{\ensuremath{\kk_l}}
\newcommand{\Latent}{\ensuremath{\ell}}
\newcommand{\alphac}{\ensuremath{\alpha}}
\newcommand{\alphaRef}{\ensuremath{\tilde \alphac}}
\newcommand{\weight}{\ensuremath{\lambda}}

\providecommand{\temp}{} 
\renewcommand{\temp}{\ensuremath{\varTheta}}
\newcommand{\tempRef}{\ensuremath{\tilde \temp}}
\newcommand{\tempDiff}{\ensuremath{\temp_\Delta}}
\newcommand{\tempDiffbar}{\ensuremath{\bar{\temp}_\Delta}}
\newcommand{\tempDiffdot}{\ensuremath{\dot{\temp}_\Delta}}
\newcommand{\tempDiffGamma}{\ensuremath{\temp_{\Delta,\Gamma}}}
\newcommand{\tempGamma}{\ensuremath{\temp_{\Gamma}}}
\newcommand{\Vint}{\ensuremath{\V_{\text{int}}}}
\newcommand{\V}{\ensuremath{\Upsilon}}
\newcommand{\VRef}{\ensuremath{\tilde{\V}}}
\newcommand{\VDiff}{\ensuremath{\V_\Delta}}
\newcommand{\ofGradTemp}{\ensuremath{(\nabla \temp\time)}}

\newcommand{\ofGradTempFEM}{\ensuremath{(\nabla \temp\hFEM\time)}}
\newcommand{\ofGradTempRefFEM}{\ensuremath{(\nabla \tempRef\hFEM\time)}}

\newcommand{\ofTemp}{\ensuremath{(\temp\time)}}
\newcommand{\Vtext}{\ensuremath{\V\ofGradTemp}}
\newcommand{\control}{\ensuremath{u}}

\newcommand{\controlDiff}{\ensuremath{\control_\Delta}}

\newcommand{\controlf}{\ensuremath{\control_\cK}}
\newcommand{\intnormal}{\ensuremath{\normal_{\text{int}}}}
\newcommand{\normal}{\ensuremath{\boldsymbol{n}}}
\newcommand{\trialf}{\ensuremath{\textsf{v}\hFEM}}
\newcommand{\trialfV}{\ensuremath{\hat{\textsf{v}}\hFEM}}
\newcommand{\pert}{\ensuremath{\varphi}}
\newcommand{\controlx}[1][]{%
\ifthenelse{\isempty{#1}}{\control}{\ensuremath{\control_{#1}}}}%
\newcommand{\controltilde}[1][]{%
\ifthenelse{\isempty{#1}}{\control}{\ensuremath{\tilde{\control}_{#1}}}}%
\newcommand{\pertx}[1][]{%
\ifthenelse{\isempty{#1}}{\pert}{\ensuremath{\pert_{#1}}}}%

\newcommand{\MT}{\ensuremath{ M_{\temp}}}
\newcommand{\ATT}{\ensuremath{ A_{\temp \temp}}}
\newcommand{\ATV}{\ensuremath{ A_{\temp \V}}}
\newcommand{\AVV}{\ensuremath{ A_{\V \V}}}
\newcommand{\AVT}{\ensuremath{ A_{\V \temp}}}
\newcommand{\BT}{\ensuremath{B_\temp}}
\newcommand{\CT}{\ensuremath{C_\temp}}
\newcommand{\ATTtilde}{\ensuremath{ \tilde{A}_{\temp \temp}}}
\newcommand{\ATVtilde}{\ensuremath{ \tilde{A}_{\temp \V}}}
\newcommand{\AVTtilde}{\ensuremath{ \tilde{A}_{\V \temp}}}
\newcommand{\MTtilde}{\ensuremath{ \tilde{M}_{\temp}}}

\newcommand{\ATTGamma}{\ensuremath{ A_{\tempGamma \tempGamma}}}

\newcommand{\AVTGamma}{\ensuremath{ A_{\V \tempGamma}}}

\newcommand{\cM}{\ensuremath{\mathcal{M}}}
\newcommand{\cA}{\ensuremath{\mathcal{A}}}
\newcommand{\cB}{\ensuremath{\mathcal{B}}}
\newcommand{\cC}{\ensuremath{\mathcal{C}}}
\newcommand{\cK}{\ensuremath{\mathcal{K}}}
\newcommand{\cS}{\ensuremath{\mathcal{S}}}

\newcommand{\bX}{\ensuremath{\mathbf{X}}}
\newcommand{\cCint}{\ensuremath{\mathcal{C_{\text{int}}}}}
\newcommand{\cBorig}{\ensuremath{\hat{\cB}}}

\newcommand{\cost}{\ensuremath{\mathcal{J}}}
\newcommand{\by}{\ensuremath{\mathbf{y}}}
\newcommand{\bu}{\ensuremath{\mathbf{u}}}
\newcommand{\state}{\ensuremath{x}}

\newcommand{\xx}{\ensuremath{\state_1}}
\newcommand{\xy}{\ensuremath{\state_2}}
\newcommand{\stateDiff}{\ensuremath{\state_\Delta}}

\newcommand{\stateDiffdot}{\ensuremath{\dot{\state}_\Delta}}

\renewcommand{\output}{\ensuremath{y}}

\newcommand{\nstate}{\ensuremath{n}}
\newcommand{\ninput}{\ensuremath{m}}
\newcommand{\noutput}{\ensuremath{r}}

\newcommand{\nL}{\ensuremath{s}}

\newcommand{\ord}{\ensuremath{\wp}}
\newcommand{\nord}{\ensuremath{n_{\ord}}}

\newcommand{\ntimestep}{\ensuremath{n_t}}
\newcommand{\kiter}{\ensuremath{k}}
\newcommand{\kiterind}{\ensuremath{_\kiter}}

\newcommand{\kiterindbwd}{\ensuremath{_{\ntimestep-\kiter}}}
\newcommand{\kiterindT}{\ensuremath{^{\phantom{\transpose}}\kiterind}}
\newcommand{\kiterindTbwd}{\ensuremath{^{\phantom{\transpose}}\kiterindbwd}}

\newcommand{\ttime}{\ensuremath{t}}
\renewcommand{\time}{\ensuremath{(\ttime)}}
\newcommand{\Endtime}{\ensuremath{\ttime_\text{end}}}

\newcommand{\tk}{\ensuremath{\ttime\kiterind}}

\newcommand{\tkh}{\ensuremath{\hat{\ttime}\kiterindT}}

\newcommand{\timek}{\ensuremath{(\tk)}}
\newcommand{\timekh}{\ensuremath{(\tkh)}}

\newcommand{\diff}[1]{\ensuremath{\mathrm{d}#1}}
\newcommand{\timediff}{\ensuremath{\frac{\mathrm{d}}{\mathrm{d}t}}}

\newcommand{\mess}{\textsf{M.E.S.S.}}
\newcommand{\mmess}{\mbox{\textsf{M-}\mess{}}}
\newcommand{\transpose}{\ensuremath{^{\textsf{T}}}}

\newcommand{\discreteset}{\ensuremath{\mathscr{T}}}

\newcommand{\timestepsbwd}{\ensuremath{\discreteset_{\text{bwd}}}}

\newcommand{\timestepsfwdref}{\ensuremath{\discreteset_{\text{fwd}}^{\text{ref}}}}
\newcommand{\meshcells}{\ensuremath{\mathscr{Q}{\hFEM}}}
\newcommand{\meshcell}{\ensuremath{Q}}
\newcommand{\timestep}{\ensuremath{\tau}}
\newcommand{\timestepk}{\ensuremath{\timestep\kiterind}}

\newcommand{\hFEM}{\ensuremath{^h}}

\newlength\figureheight
\newlength\figurewidth
\newcounter{mymac@matlab}
  \newcommand{\matlab}{MATLAB%
   \ifnum\value{mymac@matlab}<1%
   \textsuperscript{\textregistered}%
   \setcounter{mymac@matlab}{1}%
   \fi%
  }
\newcommand{\fenics}{\textsf{FEniCS}}
\newenvironment{blockmatrix}{%
  \left[\vphantom{n}%
  \vcenter\bgroup\hbox\bgroup
  \tikzpicture[
    x=1\baselineskip,
    y=1\baselineskip,
    remember picture,
  ]%
}{%
  \endtikzpicture
  \egroup
  \egroup
  \vphantom{n}\right]%
} 

\newcommand*{\mblock}[1][white]{%
  \mblockaux{#1}%
}
\def\mblockaux#1(#2,#3)#4(#5,#6){%
  \draw[fill={#1},draw=white,line width=1pt]
  let \p1=(#2,#3),
      \p2=(#5,#6),
      \p3=(#2+#5,#3+#6),
      \p4=(#2+#5/2,#3+#6/2)
  in
    (\p1) rectangle (\p3)
    (\p4) node {$#4$}
  ;%
}

\def\lblockaux#1(#2,#3)#4(#5,#6)(#7){%
  \draw[fill={#1},draw=white,line width=1pt]
  let \p1=(#2,#3),
      \p2=(#5,#6),
      \p3=(#2+#5,#3+#6),
      \p4=(#2+#5/2,#3+#6/2)
  in
    (\p1) rectangle (\p3)
    (\p4) node(#7) {$#4$}
  ;%
}
\newcommand\restr[2]{{
  \left.\kern-\nulldelimiterspace 
  #1 
  \vphantom{\big|} 
  \right|_{#2} 
  }}
\makeatletter 
\def\ifEmpty#1{\def\@temp{#1}\ifx\@temp\@empty} 
\makeatother
\newcommand{\nrm}[2][]{\ensuremath{\left\lVert #2 \right\rVert\ifEmpty{#1}\else_{#1}\fi}}

\definecolor{uncontroled}{RGB}{217,95,2}
\definecolor{desiredpre}{RGB}{102,166,30}
\colorlet{desired}{desiredpre!20!black}
\definecolor{controled}{RGB}{117,112,179}
\definecolor{extracolor1}{RGB}{27,158,119}
\definecolor{extracolor2pre}{RGB}{230,171,2}
\colorlet{extracolor2}{extracolor2pre!95!black}
\definecolor{extracolor3}{RGB}{231,41,138}
\colorlet{inputcolor}{orange}
\definecolor{outputcolor}{HTML}{33a5c3}
\colorlet{pertcolor}{red!70!black}
\colorlet{bdf1color}{uncontroled}
\colorlet{bdf2color}{desired}
\colorlet{bdf3color}{controled}
\colorlet{bdf4color}{extracolor1}
\colorlet{split1color}{extracolor2}
\colorlet{split2color}{extracolor3}
\colorlet{nofeedbackcolor}{pertcolor}
\colorlet{fterr8color}{blue}
\colorlet{fterr9color}{green}
\colorlet{fterr10color}{orange}
\colorlet{ftu2color}{red}
\colorlet{ftdtu2color}{cyan}
\tikzset{cross/.style={cross out, draw, 
         minimum size=2*(#1-\pgflinewidth), 
         inner sep=0pt, outer sep=0pt}}
\DeclareMathOperator*{\argmax}{argmax}
\crefname{line}{Line}{Lines}
\crefalias{AlgoLine}{line}%
\Crefname{algocf}{Algorithm}{Algorithms}
\makeatletter
\let\cref@old@stepcounter\stepcounter
\def\stepcounter#1{%
  \cref@old@stepcounter{#1}%
  \cref@constructprefix{#1}{\cref@result}%
  \@ifundefined{cref@#1@alias}%
    {\def\@tempa{#1}}%
    {\def\@tempa{\csname cref@#1@alias\endcsname}}%
  \protected@edef\cref@currentlabel{%
    [\@tempa][\arabic{#1}][\cref@result]%
    \csname p@#1\endcsname\csname the#1\endcsname}}
\makeatother

\begin{document}
\title{Riccati-feedback Control of a Two-dimensional Two-phase Stefan Problem}

\author[$\ast$]{Bj\"orn Baran}
\affil[$\ast$]{Max Planck Institute for Dynamics of Complex Technical Systems, Sandtorstr. 1, 39106 Magdeburg.\authorcr
  \email{baran@mpi-magdeburg.mpg.de}, \orcid{0000-0001-6570-3653}
}
\author[$\dagger$]{Peter Benner}
\affil[$\dagger$]{Max Planck Institute for Dynamics of Complex Technical Systems, Sandtorstr. 1, 39106 Magdeburg.\authorcr
  \email{benner@mpi-magdeburg.mpg.de}, \orcid{0000-0003-3362-4103}
}
\author[$\ddagger$]{Jens Saak}
\affil[$\ddagger$]{Max Planck Institute for Dynamics of Complex Technical Systems, Sandtorstr. 1, 39106 Magdeburg.\authorcr
  \email{saak@mpi-magdeburg.mpg.de}, \orcid{0000-0001-5567-9637}
}

\shortdate{}
  


%
\shorttitle{Riccati-Feedback Control of a Stefan Problem}
\shortauthor{B. Baran, P. Benner, J. Saak}

\abstract{%
We discuss the feedback control problem for a two-dimensional two-phase Stefan problem. 
In our approach, we use
a sharp interface representation in combination with mesh-movement to track the interface position.
To attain a feedback control, we apply the linear-quadratic regulator
approach to a suitable linearization of the problem.
We address details regarding the discretization and the interface representation therein.
Further, we document the matrix assembly to generate a non-autonomous generalized differential Riccati equation.
To numerically solve the Riccati equation, we use low-rank factored and
matrix-valued versions of the non-autonomous backward differentiation formulas,
which incorporate implicit index reduction techniques.
For the numerical simulation of the feedback controlled Stefan problem, we use a time-adaptive fractional-step-theta scheme.

We provide the implementations for the developed methods and test these in several numerical experiments.
With these experiments we show that our feedback control approach is applicable to the Stefan control problem
and makes this large-scale problem computable.
Also, we discuss the influence of several
controller design parameters, such as the choice of inputs and outputs.
}

\novelty{%
We propose a linear-quadratic regulator approach for the feedback stabilization
of a two-dimensional two-phase Stefan problem, where the control target is to 
steer the interface position.
This closed-loop control problem goes beyond existing open
and closed-loop control approaches for one-dimensional or one-phase Stefan problems.
Further, this is the first time the linear-quadratic regulator approach is
applied to this type of problem with a moving interface or inner boundary.

Our approach can handle the non-linearities and differential-algebraic structures
induced by the Stefan problem as well as time-dependent matrices that are present
in a non-autonomous differential Riccati equation.
We use a new non-autonomous backward differentiation formula method 
to numerically solve this Riccati equation and compute feedback controls,
which successfully stabilize the interface position.
}

\maketitle

\section{Introduction}
The solidification and melting of materials is an active and intensively studied field with numerous applications.
This phase-change problem can be modeled by a non-linear PDE and is often called Stefan problem after J.~Stefan
who describes it in his works~\cite{Ste89a,Ste89b,Ste90}.
In a certain domain, the temperature of the material is either below, above, or equal to
the specific melting temperature of the material.
Accordingly, the domain is split into a solid and a liquid phase, which are separated by an interface or inner boundary.
J.~Stefan formulates what is now called Stefan condition in~\cite{Ste89a}, which couples the time-derivative of the interface position
with the jump of the temperature gradient along the interface.
On the other hand, J.~Stefan was not the first to consider this type of problem.
G.~Lam{\'e} and B.~P.~Clapeyron were concerned with this in an earlier work as well~\cite{LamC31},
such that it is also called Lam{\'e}--Clapeyron problem.

Several books address the Stefan problem, e.g.,~\cite{Rub71,NieCM11,Gup18,KogK20a}.
An extensive historical survey of the Stefan problem can be found in the book by L.~I.~Ruben\v{s}te\u{\i}n~\cite[Introduction: \S 1]{Rub71}.
Early works on the Stefan problem usually consider the one-dimensional case, while higher dimensional cases where only studied decades after the original publication.

The aim of this manuscript is to compute and apply feedback control to the two-dimensional two-phase Stefan problem.
While this problem has been studied in combination with open-loop controls several times, see 
e.g.~\cite{Zie08, Ber10, AntNS14, AntNS15, BarBHetal18a} and the reference therein, only
recently, also closed-loop, i.e. feedback control, for the Stefan problem has been discussed in~\cite{KogK19,KogK20,KogK20a}.
However, these works address only the one-dimensional case.
On the one hand, the novelty of our work lies in the consideration of feedback control for the two-phase Stefan problem in two spatial dimensions 
in contrast to the one-phase or one-dimensional Stefan problem.
On the other hand, we apply the linear-quadratic regulator (LQR) approach to the Stefan problem, which goes beyond the types of problems
that have been studied in connection with this approach.

In particular, the application of LQR requires the treatment of the non-linearities, the differential-algebraic nature, 
and the time-varying character of the Stefan problem.
In this manuscript, we address the details on how to transform it into a linear ordinary differential equation,
apply the LQR approach to compute feedback controls, and then use these feedback controls
to stabilize the interface position in the non-linear differential-algebraic problem.
In order to do this, we linearize the Stefan problem around a reference trajectory and assemble a generalized differential Riccati equation (DRE).
For the derivation of the LQR problem and the resulting DRE we refer to~\cite{Rei72, BitG91}.
An extensive study on DREs and their numerous applications can be found in~\cite{AboFIetal03}.
Since the Stefan condition is an algebraic equation, which is coupled to the Stefan problem, 
the DRE we are considering results from a differential-algebraic equation (DAE).
Details on generalized DREs can be found in~\cite{KunM90} and for DREs resulting from DAEs, see~\cite{KunM89}.

To solve the Stefan problem numerically and assemble the DRE, we discretize in space using the
finite element method (FEM) by applying the software \fenics~\cite{dolfin}.
This results in a large-scale matrix-valued DRE with time-dependent coefficients, denoted as a non-autonomous DRE.\
While we observe in numerical experiments that, usually, its solution can be well
approximated by a low-rank factorization, this is proven theoretically
only for the autonomous DRE in~\cite{Sti18}.
Well known methods, which use the low-rank structure of the numerical solution, 
are splitting schemes~\cite{Sti15,Sti15a,Sti18a,OstPW18,MenOPetal18},
Rosenbrock and Peer methods~\cite{Men12,LanMS15,Lan17,BenL18}
as well as the backward differentiation formulas (BDF)~\cite{BenM04, Men12, LanMS15, BenM18}.
Krylov subspace methods~\cite{BehBH21,KosM20,KirSi19,GueHJetal18} and exponential integrators~\cite{LiZL20}
for DREs have been developed recently as well.
In~\cite{Men12, BenM18} Rosenbrock and BDF methods and in~\cite{LanMS15,Lan17} also Peer methods are studied for a non-autonomous DRE
where the mass matrix is constant.
In extension of this, splitting schemes and BDF methods are developed in~\cite{BarBSS21a} for non-autonomous DREs with a time-varying mass matrix.
In contrast to the BDF methods, the splitting schemes require that 
the coefficients can be decomposed into a time-dependent scalar function times a constant matrix.
The non-autonomous DREs we consider go beyond this case.
We use mesh movement techniques to track the interface, and, as a consequence, the single matrix entries change very differently, possibly for all matrices.
Thus, the BDF methods are
the most promising method for the non-autonomous DRE resulting from the Stefan problem and we apply the non-autonomous BDF method from~\cite{BarBSS21a}.

\paragraph{Structure of the Manuscript}
In this manuscript, we derive and apply several feedback controls to the two-dimensional two-phase Stefan problem.
For this, we state the equations that model the Stefan problem and the mesh movement that we use to track the interface position (\cref{sec:stefan}).
Then, we linearize the resulting system and semi-discretize it in space (\cref{sec:discretization}).
With the resulting time-dependent matrices we formulate a non-autonomous DRE resulting from the LQR approach (\cref{sec:LQR}).
We solve this DRE with the non-autonomous BDF method and use the resulting feedback gain matrices to compute a feedback control (also \cref{sec:LQR}).
We apply the described methods in several numerical experiments for different parameter settings and specify where all our codes
and data are available (\cref{sec:experiments}).

\paragraph{Notation}
In most equations, we omit the time-dependence \time{}, the spatial dependence $(\state)$, or the combination of both $(\ttime, \state)$.
This is supposed to improve the readability of the equations.

\section{Two-dimensional Two-phase Stefan Problem}\label{sec:stefan}
In this section, we define the equations describing the Stefan problem, its boundary conditions and initial values. 
This involves equations characterizing the temperature and interface movement.
Our goal is to
use a Riccati feedback control approach 
to control the interface position.
For this,
we choose a sharp interface representation in the formulation of the Stefan problem,
opposed to, e.g., the level set representation of the interface in \cite{Ber10}.
As in~\cite{BaePS10,BaePS13}, we extend the interface movement to the whole domain.
This ensures the mesh regularity for the semi-discretized Stefan problem that is described in \cref{sec:discretization}.

At each \(\ttime \in [0, \Endtime]\), the domain is $\domain\time \subset
\mathbb{R}^2$. 
One instance is illustrated in \cref{fig:domain}. 
We split the domain \domain\time{} into the two regions corresponding to the two phases.
These are the region where the material is in its solid phase \domains\time{} and accordingly, the region \domainl\time{} related to the liquid phase.
The two phases are separated by the interface \Gammaint\time.
This inner phase-boundary can move such that its position is time-dependent.
Thus, also the two phases \domains\time{} and \domainl\time{} are time-dependent 
and, as a consequence, so is the whole domain \domain\time{} and its boundary regions.
The boundary of \domain\time{} is separated into \Gammain\time{}, \Gammacool\time{} and \GammaN\time{} 
as depicted in \cref{fig:domain}.
Note that the outer shape of \domain\time{} is constant for the realization chosen in this manuscript.
Thus, the time-dependence of \domain\time{} is not absolutely necessary even though
its sub-domains are time-dependent.
However, we keep the time-dependence in our notation in order to not restrict our methods 
to this case.

  \tikzexternalenable%
  \tikzsetnextfilename{domain}%
  \filemodCmp{tikz/domain.tikz}{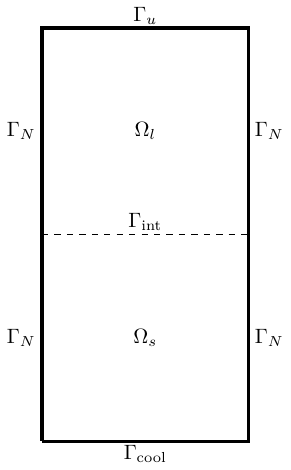}%
  {\tikzset{external/remake next}}{}%
  \input{tikz/domain.tikz}%
  \tikzexternaldisable%

We follow the definition of the Stefan problem from~\cite{BarBHetal18a}. However, here we use
it in a more compact form, i.e.\ omitting the couplings with the Navier-Stokes equations and the interface graph formulation.
While the Navier-Stokes equations alone add additional algebraic constraints, making the DAE harder to classify,
both the Navier-Stokes equations and the interface graph formulation add more nonlinearities to the problem.
In order to develop the general numerical strategy, 
we first want to study the feedback control problem for this simplified setting without these couplings,
before delving into the additional technical challenges of the full problem formulation.

We denote the temperature as \temp\time{} and model it with the partial
differential \cref{eq:heat}:
\begin{subequations}\label{eq:heat}\begin{align}
    \label{eq:heat1}
    \dot\temp - \V \cdot \nabla \temp - \alphac \Delta \temp	&= 0,  				& &\text{on}\ (0, \Endtime] \times \domain, \\
    \label{eq:heat2}
    \partial_{\normal}\temp						&= \control,			& &\text{on}\ (0, \Endtime] \times \Gammain, \\
    \label{eq:heat3}
    \temp								&= \Tcool,			& &\text{on}\ (0, \Endtime] \times \Gammacool, \\
    \label{eq:heat4}
    \temp								&= \Tmelt,			& &\text{on}\ (0, \Endtime] \times \Gammaint, \\
    \label{eq:heat5}
    \partial_{\normal} \temp						&= 0,				& &\text{on}\ (0, \Endtime] \times \GammaN, \\
    \label{eq:heat6}
    \temp(0)								&= \Tzero,			& &\text{on}\ \domain.
\end{align}\end{subequations}

In \cref{eq:heat2}, we apply the control \control\time{} as a Neumann condition on the control boundary \Gammain\time.
The \cref{eq:heat3,eq:heat4} describe the Dirichlet conditions on the cooling boundary \Gammacool\time{} and the interface \Gammaint\time{} 
with the constants \Tcool{} and \Tmelt{}, respectively. 
\cref{eq:heat6} presents the initial condition with the initial temperature distribution \Tzero.
The heat conductivities in the solid phase \ks{} and in the liquid phase \kl{} are collected in \alphac\ofTemp{}:
\begin{equation*}
    \alphac = \begin{cases}
		      \ks, & \text{on}\ \domains,\\
		      \kl, & \text{on}\ \domainl.\\
               \end{cases}
\end{equation*}
In \cref{eq:heat1}, the temperature is coupled with the extended interface movement \Vtext{}. 
For every $\ttime \in (0, \Endtime]$, we model \Vtext{} with a system of algebraic equations,
in the sense that they do not contain any time-derivatives:

\begin{subequations}\label{eq:V}\begin{align}
    \label{eq:V1}
    \Delta \V		&= 0, 									& &\text{on}\ \domain, \\
    \label{eq:V2}
    \V 			&= \Big(\frac{1}{\Latent}[\kk(\nabla \temp\time)]_l^s\Big) \cdot \intnormal,	& &\text{on}\ \Gammaint,\\
    \label{eq:V3}
    \V			&= 0, 									& &\text{on}\ \Gammacool \cup \Gammain, \\
    \label{eq:V4}
    \V \cdot \normal	&= 0, 									& &\text{on}\ \GammaN, \\
    \label{eq:V5}
    \V(0)			&= 0, 									& &\text{on}\ \domain, 
\end{align}\end{subequations}

On the interface, $\restr{\V}{\Gammaint}\ofGradTemp = \Vint\ofGradTemp$ is the interface movement in normal direction,
where \intnormal\time{} is the unit normal vector pointing from \domains\time{} to \domainl\time.
\Vint\ofGradTemp{} is coupled to \temp\time{} through the Stefan condition~\eqref{eq:V2}.
Here, \Latent{} is the latent heat constant and
\begin{equation}\label{stefan}
   [k(\nabla \temp)]_l^s = \ks\partial_{\intnormal}\restr{\temp}{\domains} - \kl\partial_{-\intnormal}\restr{\temp}{\domainl}
\end{equation}
is the jump of the temperature gradient along \Gammaint\time.
\cref{eq:V3,eq:V4} ensure that the outer boundaries of \domain\time{} do not move
such that the outer shape of the domain does not change.

We use \Vtext{} in \cref{sec:discretization} for the mesh movement.
The Stefan problem, which is described by the system of \cref{eq:heat,eq:V},
is a non-linear system of DAEs on a
time-varying domain.
Due to the coupling of \temp\time{} and \Vint\ofGradTemp{} in~\eqref{eq:V2}
and the two temperature-dependent phases in the definition of \alphac\ofTemp{},
both, \Vint\ofGradTemp{} and \alphac\ofTemp{}, depend on \temp\time{}.
Thus, the terms $\Vtext \cdot \nabla \temp\time$ and $\alphac\ofTemp \Delta \temp\time$ in \cref{eq:heat1} are nonlinear.

\section{Discretization and Linearization}\label{sec:discretization}
In order to apply the LQR approach in \cref{sec:LQR}, 
we need to formulate the Stefan problem in a standard state-space format,
which is linear and semi-discretized in space.
Thus, we describe how to transform the coupled DAE system of \cref{eq:heat,eq:V}
into
\begin{equation}\label{eq:sss}\begin{aligned}
  \cM\dot\state\hFEM 	&= \cA\state\hFEM + \cBorig\control\hFEM, \\
  \output\hFEM 			&= \cC\state\hFEM,
\end{aligned}\end{equation} 
in this section.
To generate the square matrices $\cA\time, \cM\time \in \mathbb{R}^{\nstate \times \nstate}$, 
the input matrix $\cBorig\time \in \mathbb{R}^{\nstate \times \ninput}$, and
output matrix $\cC\time \in \mathbb{R}^{\noutput \times \nstate}$,
we spatially discretize and linearize the Stefan problem.
Further, we take special care on how the boundary conditions \labelcref{eq:heat4,eq:V2} 
are treated in the definition of the matrices for \cref{eq:sss}
since they are of particular importance for our feedback control problem.


For the spatial discretization, we use the finite element method (FEM) 
on a mesh of triangular cells $\meshcells\time = \{\meshcell\time\}$ that changes over time, driven by the movement of the interface.
The interface \Gammaint\time{} itself is represented explicitly and sharply
through facets that are aligned with \Gammaint\time.
In order to track \Gammaint\time{} with the mesh, we use the semi-discrete
extended interface movement \V\hFEM\ofGradTempFEM{} to adapt the mesh inside
the whole domain, as described in~\cite{BaePS10,BaePS13}.
In this way, we prevent mesh tangling and too strong deformations.
See~\cite{BarBHetal18a} for a detailed description of our implementation.

In order to define the linearization of the Stefan problem, we require a trajectory generated
by applying an open-loop control approach to the nonlinear problem. For this, we take the open-loop control approach from~\cite{BarBHetal18a}.
This (desired) reference trajectory contains the semi-discrete reference solutions \tempRef\hFEM\time, \VRef\hFEM\ofGradTempRefFEM, \GammaintReft, and
the heat conductivities \alphaRef\time, depending on the reference trajectory.

Combining these two techniques, we derive the linearized, semi-discrete version of \cref{eq:heat1}
\begin{equation}\label{eq:heatRef}
  \dot\temp\hFEM - \V\hFEM \cdot \nabla \tempRef\hFEM - \alphaRef \Delta \temp\hFEM = 0, \qquad\text{on}\ (0, \Endtime] \times \domain, 
\end{equation}
by using the semi-discrete states ($\temp\hFEM\time, \V\hFEM\ofGradTempFEM$) and, in particular, 
by replacing \temp\hFEM\time{} in the convection term by \tempRef\hFEM\time{} and \alphac\time{}, which
depends on \temp\hFEM\time{}, by the reference heat conductivity \alphaRef\time{}.

In order to get the matrices for \cref{eq:sss}, we pose the semi-discrete variational formulations 
of \cref{eq:heatRef} together with its boundary conditions and \cref{eq:V}:
\begin{align*}
    0 &= \int\limits_{\domain} \dot\temp\hFEM \cdot \trialf \mathrm{d}x
	- \int\limits_{\domain} \V\hFEM \cdot \nabla \tempRef\hFEM \cdot \trialf \mathrm{d}x
	+ \int\limits_{\domain} \alphac \nabla\temp\hFEM\cdot\nabla\trialf \mathrm{d}x
	- \int\limits_{\Gammain} \kl \control\hFEM \cdot \trialf \mathrm{d}s,\\
    0 &= - \int\limits_{\domain} \nabla\V\hFEM\cdot\nabla\trialfV  \mathrm{d}x
	  + \int\limits_{\Gammaint} \V\hFEM \cdot \trialfV \mathrm{d}s
	  - \int\limits_{\Gammaint} \frac{1}{\Latent}[\kk(\nabla \temp\hFEM)]_l^s \cdot \intnormal \cdot \trialfV \mathrm{d}s.
\end{align*}
We denote the semi-discrete test functions with \trialf\time{} and \trialfV\time{}
and reformulate these variational formulations into a matrix-based form:
\begin{align*}
    \dot\temp\hFEM		&= \ATT\temp\hFEM + \ATV\V\hFEM + \BT\control\hFEM, 	& &\text{on}\ (0, \Endtime] \times \domain, \\
    0					&= \AVT\temp\hFEM + \AVV\V\hFEM, 						& &\text{on}\ (0, \Endtime] \times \domain.
\end{align*}
Here, the coefficient matrices are defined via the inner products
\begin{equation}\label{eq:blockmatrix}\begin{aligned}
    \langle \MT\temp\hFEM, \trialf \rangle 	&= \int\limits_{\domain} \temp\hFEM \cdot \trialf\, \mathrm{d}x, \\
    \langle \ATT\temp\hFEM, \trialf \rangle 	&= - \int\limits_{\domain} \alphac \nabla \temp\hFEM \cdot \nabla \trialf\, \mathrm{d}x, \\
    \langle\ATV\V\hFEM, \trialf \rangle	&= \int\limits_{\domain} \V\hFEM \cdot \nabla \tempRef\hFEM \cdot \trialf\, \mathrm{d}x, \\
    \langle \AVV\V\hFEM, \trialfV \rangle	&= \int\limits_{\domain} \nabla\V\hFEM\cdot\nabla\trialfV \mathrm{d}x
						- \int\limits_{\Gammaint} \V\hFEM \cdot \trialfV \mathrm{d}s \\
    \langle \AVT\temp\hFEM, \trialfV \rangle	&= \int\limits_{\Gammaint} \frac{1}{\Latent}[\kk(\nabla \temp\hFEM)]_l^s \cdot \intnormal \cdot \trialfV \mathrm{d}s, \\
    \langle \BT\control\hFEM, \trialf \rangle 	&= \int\limits_{\Gammain} \kl \control\hFEM \cdot \trialf \mathrm{d}s.
\end{aligned}\end{equation} 
With these definitions, we can formulate the semi-discrete linearized Stefan problem in the format of \cref{eq:sss}:
\begin{equation}\label{eq:blocksss}\begin{aligned} 
\begin{blockmatrix}
  \mblock(0,1)\MT\vphantom{\Big(}(1.6,1)	\mblock(1.6,1)0(1.6,1)
  \mblock(0,0)0\vphantom{\Big(}(1.6,1)		\mblock(1.6,0)0(1.6,1)
\end{blockmatrix}
\timediff \begin{blockmatrix}
	  \mblock(0,2)\temp\hFEM\vphantom{\Big(}(1.2,1)
	  \mblock(0,1)\V\hFEM\vphantom{\Big(}(1.2,1)
	\end{blockmatrix}
  &= \begin{blockmatrix}
  \mblock(0,1)\ATT\vphantom{\Big(}(2,1)		\mblock(2,1){\ATV}(2,1)
  \mblock(0,0)\AVT\vphantom{\Big(}(2,1)		\mblock(2,0)\AVV(2,1)
\end{blockmatrix}
\begin{blockmatrix}
	  \mblock(0,2)\temp\hFEM\vphantom{\Big(}(1.2,1)
	  \mblock(0,1)\V\hFEM\vphantom{\Big(}(1.2,1)
\end{blockmatrix}
+ 
\begin{blockmatrix}
	  \mblock(0,2)\BT\vphantom{\Big(}(1.2,1)
	  \mblock(0,1)0\vphantom{\Big(}(1.2,1)
\end{blockmatrix} 
\bu\hFEM, \\ \\
\by\hFEM\ &=
\begin{blockmatrix}
  \mblock(1,4)\CT(1,1)	\mblock(2,4)0(1,1)
\end{blockmatrix}
\begin{blockmatrix}
	  \mblock(0,1)\temp\hFEM\vphantom{\Big(}(1.2,1)
	  \mblock(0,0)\V\hFEM\vphantom{\Big(}(1.2,1)
\end{blockmatrix}.
\end{aligned}\end{equation}
Several choices are
plausible for the output matrix $\CT\time$.
We will introduce some in \cref{sec:experiments}.
With the zero-blocks in \cref{eq:blocksss}, the DAE structure is clearly visible.
Since, at every time instance \ttime{}, \AVV\time{} represents the Poisson operator with a Dirichlet 
boundary part,
it is always non-singular.
Further, the mass matrix with respect to the temperature \MT\time{} is symmetric positive definite
and, in particular, always non-singular as well.
Thus, \cref{eq:blocksss} is a DAE in semi-explicit form of differential index 1 (see e.g.~\cite{KunM06}) and we can
apply the implicit index-reduction techniques from~\cite{morFreRM08} for this type of DAEs.

Next, we discuss the treatment of the boundary conditions in the FEM matrices in \cref{eq:blockmatrix,eq:blocksss}. 
This is of special importance for our objective to use these matrices to compute a 
feedback control for the Stefan problem in \cref{sec:LQR}.
The subject of attention in our feedback control problem is the position of \Gammaint\time{},
which is directly linked to the boundary conditions \labelcref{eq:heat4,eq:V2}
since they are defined on \Gammaint\time{}.
\cref{eq:heat4,eq:V2} are also important for the coupling of \cref{eq:heat,eq:V}.
This coupling is of particular importance for the feedback control problem
since the interface movement is defined in the Stefan condition \labelcref{eq:V2}
and the control is applied in \cref{eq:heat2}.
One common approach to handle Dirichlet boundary conditions, such as \cref{eq:heat4,eq:V2}, is to
remove the rows and columns corresponding to the degrees of freedom (DOFs) 
at the related boundary regions from the FEM matrices.
However, these DOFs are related to the coupling of \temp\hFEM\time{} and \V\hFEM\ofGradTempFEM, and,
thus, the coupling of the temperature and the interface movement.
By removing these DOFs, the matrices would lose important information, which is
necessary for the computation of a feedback control.
Thus, the corresponding DOFs on \Gammaint\time{} are to be present in the FEM matrices.
In order to treat the boundary conditions \labelcref{eq:heat4,eq:V2} appropriately
for our feedback control problem
we, from now on, consider the difference states
\begin{align*}
  \tempDiff\hFEM 	&= \temp\hFEM - \tempRef\hFEM, \\
  \VDiff\hFEM 		&= \V\hFEM - \VRef\hFEM.
\end{align*} 
We formulate our feedback control problem in terms of $\tempDiff\hFEM\time$ and $\VDiff\hFEM\ofGradTempFEM$.
As a result, the desired state, which we stabilize with a feedback control, is the all-zero-state.
The semi-discrete version of the Dirichlet condition in \cref{eq:heat4} in terms of the difference state $\tempDiff\hFEM\time$ is
equivalent to 
\begin{equation}\label{eq:BCequiv}
  0 = \TmeltDiff - \tempDiff\hFEM, \qquad\text{on}\ (0, \Endtime] \times \Gammaint.
\end{equation}
Further, both terms $\TmeltDiff\time = \Tmelt\time - \TmeltRef\time = 0$ and $\restr{\tempDiff\hFEM\time}{\Gammaint\time} = 0$ in 
\cref{eq:BCequiv} are constant.
Thus, also the time derivative equals zero, $\tempDiffdot\hFEM\time = 0$,
and we add this equation to \cref{eq:BCequiv}.
As a result, we can formulate the modified condition \labelcref{eq:BCDiff} to replace \cref{eq:heat4} in terms of the difference state 
for the FEM matrices, yielding
\begin{equation}\label{eq:BCDiff}
  \tempDiffdot\hFEM = \TmeltDiff - \tempDiff\hFEM, \qquad\text{on}\ (0, \Endtime] \times \Gammaint.
\end{equation}
To incorporate \cref{eq:BCDiff} into the matrices, we denote $\tempDiffGamma\hFEM\time = \restr{\tempDiff\hFEM\time}{\Gammaint\time}$ and assume we can split
\begin{align*}
\tempDiff\hFEM = 
\begin{blockmatrix}
	  \mblock(0,1)\tempDiffbar\hFEM\vphantom{\Big(}(1.2,1)
	  \mblock(0,0)\tempDiffGamma\hFEM\vphantom{\Big(}(1.2,1)
\end{blockmatrix}.
\end{align*}
In order to ensure the Dirichlet condition on \Gammaint\time{}, we add \cref{eq:BCDiff} to the related block matrices, which then read
\begin{equation*}
  \MTtilde = \begin{blockmatrix}
  \mblock(0,2)\MT\vphantom{\Big(}(1.2,1)	\mblock(1.2,2)0(1.2,1)
  \mblock(0,1)0\vphantom{\Big(}(1.2,1)		\mblock(1.2,1)I(1.2,1)
\end{blockmatrix},\ 
  \ATTtilde = \begin{blockmatrix}
  \mblock(0,2)\ATT\vphantom{\Big(}(2,1)		\mblock(2,2)\ATTGamma(2,1)
  \mblock(0,1)0\vphantom{\Big(}(2,1)		\mblock(2,1){-I}(2,1)
\end{blockmatrix},\ 
  \ATVtilde = \begin{blockmatrix}
  \mblock(0,2)\ATV\vphantom{\Big(}(2,1)
  \mblock(0,1)0\vphantom{\Big(}(2,1)
\end{blockmatrix},\ 
  \AVTtilde = \begin{blockmatrix}
  \mblock(0,2)\AVT\vphantom{\Big(}(2,1)		\mblock(2,2)\AVTGamma(2,1)
\end{blockmatrix}.
\end{equation*}
With this, the conditions in \cref{eq:heat2,eq:V2} can be incorporated explicitly into the definition of
\AVT\time{} and \BT\time.
Regarding the remaining boundary conditions, the Neumann condition in \cref{eq:heat5} does not need any extra attention, since it is incorporated
automatically. 
The Dirichlet boundary conditions in \cref{eq:heat3,eq:V3,eq:V4} can be handled by the common approach
to remove the rows and columns corresponding to the DOFs at the related boundary regions from the matrices.

\cref{eq:blocksss} with \cref{eq:BCDiff} incorporated reads
\begin{equation}\label{eq:blocksssGamma}\begin{aligned}  
\begin{blockmatrix}
  \mblock(0,1)\MTtilde\vphantom{\Big(}(1.6,1)	\mblock(1.6,1)0(1.6,1)
  \mblock(0,0)0\vphantom{\Big(}(1.6,1)		\mblock(1.6,0)0(1.6,1)
\end{blockmatrix}
\timediff \begin{blockmatrix}
	  \mblock(0,2)\tempDiff\hFEM\vphantom{\Big(}(1.2,1)
	  \mblock(0,1)\VDiff\hFEM\vphantom{\Big(}(1.2,1)
	\end{blockmatrix}
  &= \begin{blockmatrix}
  \mblock(0,1)\ATTtilde\vphantom{\Big(}(2,1)		\mblock(2,1){\ATVtilde}(2,1)
  \mblock(0,0)\AVTtilde\vphantom{\Big(}(2,1)		\mblock(2,0)\AVV(2,1)
\end{blockmatrix}
\begin{blockmatrix}
	  \mblock(0,2)\tempDiff\hFEM\vphantom{\Big(}(1.2,1)
	  \mblock(0,1)\VDiff\hFEM\vphantom{\Big(}(1.2,1)
\end{blockmatrix}
+ 
\begin{blockmatrix}
	  \mblock(0,2)\BT\vphantom{\Big(}(1.2,1)
	  \mblock(0,1)0\vphantom{\Big(}(1.2,1)
\end{blockmatrix} 
\bu\hFEM, \\ \\
\by\hFEM\ &=
\begin{blockmatrix}
  \mblock(1,4)\CT(1,1)	\mblock(2,4)0(1,1)
\end{blockmatrix}
\begin{blockmatrix}
	  \mblock(0,1)\tempDiff\hFEM\vphantom{\Big(}(1.2,1)
	  \mblock(0,0)\VDiff\hFEM\vphantom{\Big(}(1.2,1)
\end{blockmatrix}.
\end{aligned}\end{equation}
Notably, this modification preserves the DAE structure and index of \cref{eq:blocksss}.
Consequently, we follow~\cite{morFreRM08} and apply the Schur complement to remove the algebraic conditions. 
This yields an equivalent formulation, a.k.a. realization, 
of the Stefan problem on the hidden manifold \cite{KunM06},
i.e., as an ODE, with the coefficient matrices 
\begin{equation}\label{eq:StefanMats}\begin{aligned}
  \cM\time &= \MTtilde\time, \\
  \cA\time &= \ATTtilde\time - \ATVtilde\time\AVV^{-1}\time\AVTtilde\time, \\
  \cBorig\time &= \BT\time, \\
  \cC\time &= \CT\time.
\end{aligned}\end{equation}
With the matrices from \cref{eq:StefanMats}, we can transform the Stefan problem into the formulation of \cref{eq:sss}
and are able to apply the LQR approach for the computation of a feedback control.
In order to have a computationally efficient method, we keep the sparse structure of the matrices and
never compute the generally dense Schur complement for \cA\time{} explicitly.
Instead we apply \cA\time{} implicitly and use the matrices from \cref{eq:blocksssGamma}.
For details, see~\cite{morFreRM08}, or the numerical implementation in~\cite{SaaKB21-mmess-2.1}.

In \cref{sec:LQR}, we describe the related numerical methods.
For this, we further discretize the Stefan problem in time.
For the simulation forward in time, we use the reference time-steps
\begin{equation}\label{eq:timesteps}\begin{aligned}
  0 &= \ttime_0 < \ttime_1 < \ldots < \ttime_{\ntimestep - 1} < \ttime_{\ntimestep} = \Endtime, \\
  \timestepsfwdref &= \{\ttime_0, \ttime_1, \ldots,\ttime_{\ntimestep - 1}, \ttime_{\ntimestep}\}.
\end{aligned}\end{equation}
In \cref{{sec:LQR}}, we solve differential Riccati equations backwards in time.
For this we use the same time-steps $\hat \ttime_{k} = \ttime_{\ntimestep - k}$ in reversed order
\begin{equation}\label{eq:timestepsbwd}\begin{aligned}
  \Endtime &= \hat \ttime_0 > \hat \ttime_1 > \ldots > \hat \ttime_{\ntimestep - 1} > \hat \ttime_{\ntimestep} = \ttime_0, \\
  \timestepsbwd &= \{\hat \ttime_0, \hat \ttime_1, \ldots,\hat \ttime_{\ntimestep - 1}, \hat \ttime_{\ntimestep}\}.
\end{aligned}\end{equation}
We use the same time-steps because the matrices from \cref{eq:StefanMats},
which form the coefficients of the
differential Riccati equations, are
assembled during the forward simulation.
When we apply a time-adaptive method in the forward simulation, additional time-steps can be added to \timestepsfwdref{}.

\section{Non-autonomous Linear-Quadratic Regulator}\label{sec:LQR}
In this section, we formulate the Riccati-feedback approach for the Stefan problem.
We focus on the non-autonomous character of this problem, which is induced 
by the moving interface and the consequently changing sub-domains.
This results in time-dependent coefficients in \cref{eq:sss}.

In order to derive a feedback-based stabilization of the Stefan problem,
we use the LQR approach (e.g.,~\cite{Son98}). 
We use this approach because it is well studied for related types of problems, e.g., convection diffusion equations~\cite{Wei16}, and demonstrates promising performance for these.

To formulate the control problem, 
we define a quadratic cost functional $J(\output\hFEM,\controlDiff\hFEM)$, 
tracking the deviation of the output
from the desired output as well as penalizing the control costs with
a weight factor $ 0 <\weight \in \mathbb{R}$. The cost functional is thus defined as
\begin{equation}\label{eq:cost} 
  \cost(\output\hFEM,\controlDiff\hFEM) = \frac{1}{2}\int_0^{\Endtime} \nrm{\output\hFEM - \output\hFEM_{d}}^2 + \weight \nrm{\controlDiff\hFEM}^2\;\diff t. 
\end{equation}
We minimize this cost functional subject to 
the linear time-varying system~\labelcref{eq:sss}.
With the matrices from \cref{eq:StefanMats} the LQR problem reads
\begin{equation}\begin{aligned}\label{eq:LQR}
	&\min\limits_{\controlDiff\hFEM} J(\output\hFEM, \controlDiff\hFEM) \\
	\text{subject to} &\\
	&\cM\stateDiffdot\hFEM 	= \cA\stateDiff\hFEM + \cBorig\controlDiff\hFEM, \\
  	&\phantom{\cM\stateDiffdot\hFEM}\makebox[0pt][r]{\output\hFEM} 			= \cC\stateDiff\hFEM.
\end{aligned}\end{equation}
The unique solution to the LQR problem \labelcref{eq:LQR} is 
(see \cref{thrm:unique})
\begin{equation}\label{eq:u}
  \controlDiff\hFEM = -\cK \stateDiff\hFEM,
\end{equation}
where the feedback gain matrix
\begin{equation}\label{eq:K}
  \cK = \frac{1}{\weight}\cBorig\transpose\bX\cM = \frac{1}{\sqrt{\weight}}\cB\transpose\bX\cM
\end{equation}
requires the solution $\bX\time \in \mathbb{R}^{\nstate \times \nstate}$
of a differential Riccati equation (DRE).
For the Stefan problem, this is the large-scale matrix-valued non-autonomous generalized DRE
\begin{equation}\label{eq:DREpre}
      - \timediff{(\cM\transpose \bX \cM)} = \cC\transpose\cC+\cA\transpose\bX\cM+\cM\transpose \bX\cA
      -\cM\transpose\bX\cB\cB\transpose\bX\cM.
\end{equation}
All coefficients of \cref{eq:DREpre} can be non-autonomous (but we skip the $(\ttime)$-dependency for better readability).
The coefficients of the DRE, at each time instance \ttime{}, are the matrices 
$\cA\time, \cM\time$, $\dot\cM\time \in \mathbb{R}^{\nstate \times \nstate}$, $\cB\time \in \mathbb{R}^{\nstate \times \ninput}$, 
and $\cC\time \in \mathbb{R}^{\noutput \times \nstate}$ from \cref{eq:StefanMats},
where the input matrix $\cB\time = \frac{1}{\sqrt{\weight}}\cBorig\time$ is scaled.
In order to solve the DRE, the time-derivative requires special treatment due to the time-dependent mass matrix \cM\time{}.
The left-hand side in \cref{eq:DREpre} can be computed applying the chain rule,
\begin{equation*}
 	- \timediff{(\cM\transpose \bX \cM)} = - \dot\cM\transpose \bX \cM - \cM\transpose \dot\bX \cM - \cM\transpose \bX \dot\cM.
\end{equation*}
We subtract the two terms containing $\dot\cM\time$ from \cref{eq:DREpre} to obtain a DRE 
with the time-derivative of \cM\time{} moved to the right-hand side:
\begin{equation}\label{eq:DRE}
      - \cM\transpose \dot\bX \cM = \cC\transpose\cC+(\dot\cM + \cA)\transpose\bX\cM+\cM\transpose \bX(\dot\cM + \cA)
      -\cM\transpose\bX\cB\cB\transpose\bX\cM.
\end{equation}

\begin{assumption}\label{ass:reg}
	The matrix pencil $\alpha\cM\time - \beta(\dot\cM\time + \cA\time)$ 
	is regular.
\end{assumption}
Note that we observed \cref{ass:reg} to hold in our numerical experiments.

\begin{theorem}\label{thrm:unique}
	If \cref{ass:reg} holds, the unique solution to the LQR problem
	\labelcref{eq:LQR} given by the control function $\controlDiff\hFEM$
	defined by \crefrange{eq:u}{eq:DRE}.
\end{theorem}
\begin{proof}
We use the result from \cite{KunM90} that the LQR problem \labelcref{eq:LQR}
can be reduced to a standard control problem under certain conditions.
These are \cite[condition (3.4) and (3.5)]{KunM90}.

With the block structure of the matrices in \cref{eq:blocksssGamma}, 
condition \cite[condition (3.5)]{KunM90} is fulfilled even without
the singular value decomposition and the resulting transformation
that is performed there. 

Together with \cref{ass:reg}, also \cite[condition (3.4)]{KunM90}
is fulfilled and \cref{eq:blocksssGamma} is already formulated in
the form of \cite[Equation (3.12)--(3.16)]{KunM90}.
Thus, the control function $\controlDiff\hFEM$
defined by \crefrange{eq:u}{eq:DRE}
is equivalent to the unique solution in
\cite[Equation (2.19)]{KunM90}.
Consequently, $\controlDiff\hFEM$ is the unique solution to
the LQR problem \labelcref{eq:LQR}.
\end{proof}

Besides the challenges arising from solving a large-scale matrix-valued DRE, 
the time-dependent matrices, especially the presence 
of $\dot\cM\time$, impose additional difficulties, 
which we address in this section.

To numerically compute the feedback gain matrix $\cK\kiterind = \cK\timek$ 
and the feedback control $\control\kiterind=\controlDiff\hFEM\timek$ for $\tk \in \timestepsbwd$,
the solution of the DRE \labelcref{eq:DRE} is required.
We use efficient low-rank methods for the computation of the numerical solution of the DRE, 
$\bX\kiterindT = \bX\timek$, since we assume that 
\bX\kiterindT{} has a low (numerical) rank motivated by~\cite{Sti18}.
For this, we use only a small number of inputs and outputs in our experiments, 
i.e.\ $\ninput,\noutput \ll \nstate$.
As introduced in~\cite{LanMS15}, which is motivated by~\cite{BenLT09}, \bX\kiterindT{} can, thus, be approximated to high accuracy by the decomposition
\begin{equation*}
  \bX\kiterindT \approx L\kiterindT D\kiterindT L\kiterind\transpose,
\end{equation*}
where the matrices $L\kiterind \in \mathbb{R}^{\nstate \times \nL}$ and $D\kiterind \in \mathbb{R}^{\nL \times \nL}$
have rank $\nL \ll \nstate$.
We use the low-rank non-autonomous backward differentiation formula (BDF) (see \cref{algo:BDF}) to solve the non-autonomous
DRE \labelcref{eq:DRE}. This method is described in greater detail
in~\cite[Section 2]{BarBSS21a}.
\IncMargin{1em}
\begin{algorithm}[t]
	\SetAlgoNoLine 
	\KwIn{$\cA\time,\cM\time,\dot\cM\time,\cB\time,\cC\time,\weight,\timestepsbwd, \ord, L_0,\ldots,L_{\ord-1}, D_0,\ldots,D_{\ord-1}$}
	\KwOut{$\cK\kiterind, \kiter = 1,\ldots,\ntimestep$}
	\caption{Non-autonomous low-rank factor BDF method of order \ord{}}
	\label{algo:BDF}
	\BlankLine
	\For{$\kiter = \ord,\ldots,\ntimestep$}{\label{line:BDF1}
	  $\cA\kiterind = \timestep\kiterind\beta(\dot\cM\timekh + \cA\timekh) - \frac{1}{2}\cM\timekh$\label{line:BDF3}\;
	  $\cM\kiterind = \cM\timekh$\label{line:BDF4}\;
	  $\cB\kiterind = \sqrt{\timestep\beta}\cB\timekh$\label{line:BDF5}\;
	  $\cC\kiterind\transpose = \left[\cC\timekh\transpose, \cM\kiterind\transpose L_{\kiter - 1}, \ldots,
	  \cM\kiterind\transpose L_{\kiter - \ord}\right]$\label{line:BDF6}\;
	  $\cS\kiterind =  
	    \begin{blockmatrix}
	      \mblock(0,3)\timestep\beta I_\noutput(3,1)	\mblock(3,3)(4,1)				\mblock(7,3)(3,1)	\mblock(10,3)(3,1)
	      \mblock(0,2)(3,1)						\mblock(3,2)-\alpha_1 D_{\kiter - 1}(4,1)	\mblock(7,2)(3,1)	\mblock(10,2)(3,1)
	      \mblock(0,1)(3,1)						\mblock(3,1)(4,1)				\mblock(7,1)\ddots(1,1)	\mblock(10,1)(3,1)
	      \mblock(0,0)(3,1)						\mblock(3,0)(4,1)				\mblock(7,0)(3,1)	\mblock(10,0)-\alpha_{\ord}D_{\kiter - \ord}(3,1)
	    \end{blockmatrix}$\label{line:BDF7}\;
	  solve ARE \labelcref{eq:ARE} for $L\kiterind$ and $D\kiterind$ \label{line:BDF8}\;
	  $\cK\kiterindTbwd = \frac{1}{\sqrt{\weight}}\cB\timekh\transpose L\kiterindT D\kiterindT L\kiterind\transpose\cM\kiterind$\label{line:BDF9}\;
	  }
\end{algorithm}
\DecMargin{1em} 
Compared to an open-loop control problem, the DRE replaces the adjoint equations in the LQR setting.
Consequently, we solve the DRE backwards in time and have the time-steps \timestepsbwd{} (see \cref{eq:timestepsbwd}) in reversed order as an input to \cref{algo:BDF}.
In each step of the BDF method, we solve an algebraic Riccati equation (ARE) 
\begin{equation}\begin{aligned}\label{eq:ARE}
      0 =\ 
      &\cC\kiterind\transpose\cS\kiterindT\cC\kiterindT
      +(\dot\cM\kiterindT + \cA\kiterindT)\transpose\bX\kiterindT\cM\kiterindT
      +\cM\kiterind\transpose \bX\kiterindT(\dot\cM\kiterindT + \cA\kiterindT)\\
      &-\cM\kiterind\transpose\bX\kiterindT\cB\kiterindT\cB\kiterind\transpose\bX\kiterindT\cM\kiterindT
\end{aligned}\end{equation}
for the low-rank solution 
$\bX\kiterindT \approx L\kiterindT D\kiterindT L\kiterind\transpose$.
The matrices in \cref{eq:ARE} are constructed in \cref{line:BDF3,line:BDF4,line:BDF5,line:BDF6,line:BDF7}.
The coefficients $\alpha_j$ and $\beta$ are taken from the literature,
e.g.~\cite{AscP98}.
The BDF method is a multistep method.
Thus, to start a BDF method of order \(\ord > 1\), the terminal values $\bX_0,\ldots,\bX_{\ord-1}$ are required with
sufficient accuracy to obtain the desired order of convergence.
These values can be computed with sufficiently small time-steps of
the order \(\ord-1\) method.
This procedure is repeated recursively for the order \(\ord-1\) method to compute terminal values for its start.
For details see~\cite[Algorithm 2]{BarBSS21a}, or the implementation in~\cite[mess\_bdf\_dre.m]{SaaKB21-mmess-2.1}.
This method requires an additional input parameter that we set to $\nord = 10$ as in~\cite{BarBSS21a}.

\cref{algo:BDF} is embedded in the open source software package \mmess{}~2.1~\cite{SaaKB21-mmess-2.1}, where it
benefits from efficient solvers for the ARE \labelcref{eq:ARE}.

To have an overview of the three steps for the Riccati-feedback stabilization of the Stefan problem,
we collect the overall procedure in \cref{algo:LQR}. 
\IncMargin{1em}
\begin{algorithm}[t]
	\SetAlgoNoLine 
	\caption{LQR for the Stefan problem}
	\label{algo:LQR}
	\BlankLine
	solve open-loop control problem to get the reference trajectory \newline
	and \cA\time, \cM\time, $\dot\cM\time$, \cB\time, \cC\time{} of the linearized problem\label{line:LQR1}\;
	solve DRE \labelcref{eq:DRE} with \cref{algo:BDF} to get $\cK\kiterind, \kiter = 1,\ldots,\ntimestep$\label{line:LQR2}\;
	apply \cK\kiterind{} in a forward simulation of the Stefan problem\label{line:LQR3}\;
\end{algorithm}
\DecMargin{1em}
The method we choose for the forward simulations of the Stefan problem in \cref{line:LQR1,line:LQR3} 
is a fractional-step-theta scheme and for solving the DRE in \cref{line:LQR2} the non-autonomous BDF methods, which we described above.
The time discretization \timestepsfwdref{}, defined in \cref{eq:timesteps}, is used in all three steps.
In \cref{line:LQR2}, it is used backwards in time (\timestepsbwd, \cref{eq:timestepsbwd}) and,
for the second forward solve in \cref{line:LQR3}, additional time steps can be added adaptively.
We use the time-adaptivity to prevent numerical instabilities 
that can occur with feedback controls that have very large variation, as demonstrated in~\cite{BarBSS21a}.
The specific time-adaptivity is tailored to include the relative change of the feedback control.
A detailed description of this time-adaptive fractional-step-theta scheme with a relative control-based indicator function
is presented in~\cite[Section 3]{BarBSS21a}.

We provide the references for the codes of the methods above in \cref{sec:code}
and showcase the behavior of \cref{algo:BDF,algo:LQR} by means of several numerical experiments.

\section{Numerical Experiments}\label{sec:experiments}
With several numerical experiments, we demonstrate that the LQR approach is 
applicable to this type of control problem, i.e., the Stefan problem, regardless
of several numerical challenges and approximations that arise
similar to other types of problems, where this approach is applied successfully.
The Stefan problem is non-linear, which we address with a linearization.
Further, the coefficients of the corresponding DRE are large-scale matrices.
Thus, we use a low-rank representation of the DRE solution to make it feasible
in regard of the computational cost and the memory requirements.
Additionally, the Stefan problem is a DAE as well as, e.g., the Navier-Stokes equations,
for which the LQR approach is applied successfully in \cite{Wei16}.
To handle the DAE structure, we use an implicit index reduction technique.

In addition,
for the Stefan problem, the matrices are time-dependent and
we approximate the time-derivative of the mass matrix by centered differences.
Consequently, we use the non-autonomous BDF method to solve the DRE and compute 
the feedback gain matrices.
Then, at intermediate time-steps that are
generated by the time-adaptive fractional-step-theta scheme,
we interpolate the feedback gain matrices.
An alternative is to interpolate the DRE solution with the corresponding low-rank factors and compute the feedback gain matrices at intermediate time-steps as in \cref{line:BDF8} of \cref{algo:BDF}.
However, this is more expensive and not necessary in our setting.
All these techniques can introduce additional approximation errors.
Still, the computed feedback controls can steer the interface back to
the desired trajectory successfully in our experiments, similar to previously studied problem types.

With several numerical experiments, we illustrate the
performance of the feedback stabilization, which is computed with \cref{algo:LQR}.
The open-loop computations in \cref{line:LQR1} are reported in~\cite{BarBHetal18a} and
the assessment of the runtime performance of different methods for the solution of the DRE in
\cref{line:LQR2} is left to~\cite{BarBSS21a}.
Instead, we focus on the behavior of the feedback stabilization in \cref{line:LQR3}.
This behavior is strongly influenced by the weight factor \weight{} in the cost
functional as well as different choices of inputs and outputs, as the experiments demonstrate. 
In order to
asses the robustness of our computed feedback controls, we are interested in uncertainties in the
cooling, i.e., we investigate different perturbations \pert\time{} to the Dirichlet
boundary condition at \Gammacool\time{}.

\begin{figure}[t]
  \centering
  \vspace{0em}\includegraphics{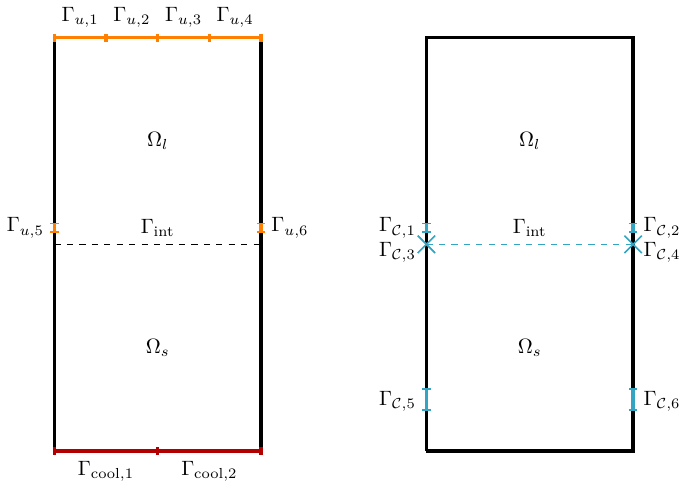}
  \caption{Input, Perturbation (left), and Output Areas (right)}%
  \label{fig:inout}
\end{figure}
We depict several boundary regions for perturbations (\Gammacool\time{}) as well as
inputs (\Gammain\time{}) and outputs (\Gammaoutx\time{}) in \cref{fig:inout}.
As for \domain\time{}, the boundary regions for inputs $\Gammainx[1]\time, \ldots, \Gammainx[4]\time$
are constant in this representation of the domain. 
However, we do not restrict ourselves to this case and, thus, keep the time-dependence in our notation.

The particular domain \domain\time{}, we choose for the experiments, is a rectangle $[0, 0.5] \times [0, 1]$
with the initial interface position at height $0.5$ as a horizontal line.
For the discretization, we choose a mesh of triangles with \(3\,899\) vertices and 401 time-steps.
This, using standard P1 elements in \fenics, results in \(3\,899\) DOFs for \temp\time{} and \(7\,798\) DOFs for \Vtext. 
After the removal of the DOFs corresponding to Dirichlet boundary conditions (see \cref{sec:discretization}),
the size of the matrices is \(\nstate{} = 11\,429\).
The model parameters are listed in \cref{tab:para}.
\begin{table}[t]
\centering
\begin{tabular}{|l|l|l|l|l|l|l|l|}
\hline
\Endtime	&	\ks	&	\kl		&	\Tcool	&	\Tmelt	&	\Latent	&	\Tzero 									&	\timestepk \\
\hline
$1$			&	$6$	&	$10$	&	$-1$	&	$0$		&	$10$	&	$\Tcool \cdot (1 - 2 \cdot \xy)$	&	 $\in \left[10^{-4}, 2.5\cdot10^{-3}\right]$\\
\hline
\end{tabular}
\caption{Stefan Problem Model Parameters}
\label{tab:para}
\end{table}

\subsection*{Code Availability}\label{sec:code}
All codes and data to reproduce the presented results are available at~\cite{LQR_Stefan_codes}.
The non-autonomous BDF method is incorporated in the software package \mmess{}~2.1~\cite{SaaKB21-mmess-2.1}

\subsection{Experiment 1}\label{sec:ex1} 
For this first experiment, we focus on the influence of weight parameters and outputs
on the performance of the feedback control.
The desired interface trajectory, which we intend to stabilize with the feedback control, 
is a flat horizontal line moving from its
initial height at \(0.5\) downward by \(0.004\). 
Since we assume the Stefan problem to be asymptotically stable,
the interface returns to the desired trajectory after a perturbation without a feedback control as well.
However, this is not achieved within the time horizon.
With the feedback control, we intend to prevent the interface deviation to a certain extend and to steer the interface back in a much shorter time after a perturbation occurs.

\begin{figure}[tp]
  \centering
  \includegraphics{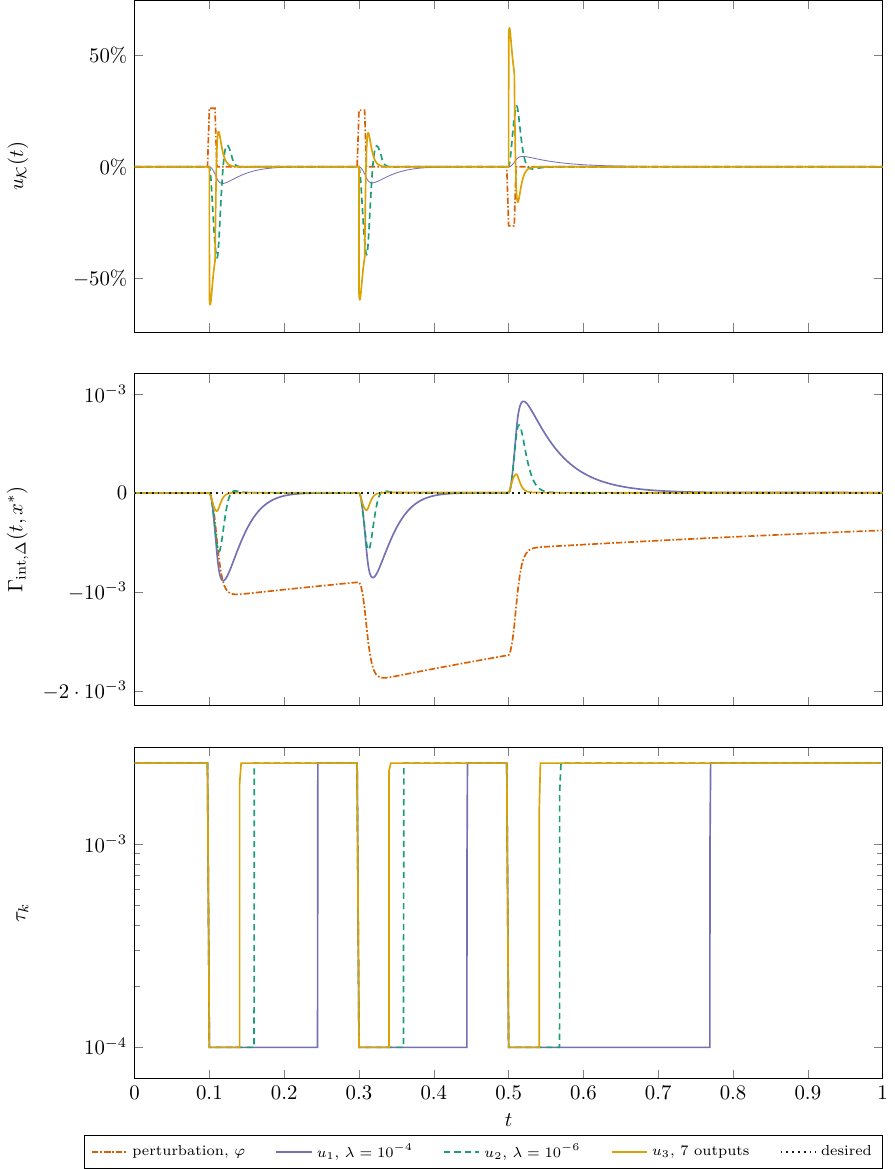}
  \caption{Perturbation and feedback (top), relative interface position (middle), 
  			and time-step sizes (bottom) for different weights and outputs
	    (\nameref{sec:ex1})}%
  \label{fig:weights}
\end{figure}

The interface deviation is caused by the perturbation \pert\time{}, which is  generated randomly in the form of three scalar values in the range of $[-\Tcool, \Tcool]$.
These are applied for a period of four time-steps of \timestepsfwdref{} at the times $0.1$, $0.3$, and $0.5$ on $\Gammacool\time =
\Gammacoolx[1]\time\cup\Gammacoolx[2]\time$:
\begin{equation*}
  \temp	= \Tcool + \pert, \qquad \text{on}\ (0, \Endtime] \times \Gammacool.
\end{equation*}
The trajectory of \pert\time{} is pictured in the top part of \cref{fig:weights}.
The scaling of the y-axis is relative to \Tcool.

Here, for the LQR problem, we use a single input \controlf\time{} with the same value on
$\Gammainx[1]\time\cup\ldots\cup\Gammainx[4]\time$, i.e.\ $\cB\time \in \mathbb{R}^{\nstate \times 1}$.
The top part of \cref{fig:weights} shows the control for three different LQR
designs, resulting from three different combinations of weights and outputs:
\begin{align*}
  &\controlx[1]: &&\weight = 10^{-4}, 		&&\text{2 outputs:}\ \Gammaoutx[3],\Gammaoutx[4],\\
  &\controlx[2]: &&\weight = 10^{-6}, 		&&\text{2 outputs:}\ \Gammaoutx[3],\Gammaoutx[4],\\
  &\controlx[3]: &&\weight = 1.6\cdot10^{-2}, 	&&\text{7 outputs:}\ \Gammaoutx[1],\ldots,\Gammaoutx[6],\Gammaint.
\end{align*}
The outputs on \Gammaoutx[1]\time{}, \Gammaoutx[2]\time{}, \Gammaoutx[5]\time{},  and \Gammaoutx[6]\time{}
measure averaged temperatures on the corresponding interval, while
\Gammaoutx[3]\time{} and \Gammaoutx[4]\time{} represent point measurements of the
temperature.
The output at \Gammaint\time{} monitors the difference of the interface movement 
to the desired movement
\begin{equation*}
  \langle \cCint , \trialf \rangle 
    = \int\limits_{\Gammaint} \Big(\frac{1}{\Latent}[\ks(\nabla \tempDiff)_s 
    - \kl(\nabla \tempDiff)_l]\Big) \cdot \intnormal \cdot \trialf \mathrm{d}s
\end{equation*}
and is defined via the jump term of the Stefan condition (\cref{eq:V2}).
This output alone is not suited for an effective control-design.
It would generate a large output, and thus an active feedback-response, only while the perturbation is actively
driving the interface away from the desired trajectory. 
However, it can not detect a difference in the position of the interface.
Consequently, it would not steer the interface back but would keep it on a ``parallel trajectory".

Thus, in the first two LQR designs, the controls \controlx[1]\time{} and \controlx[2]\time{}
are based on two outputs that measure the temperature at the desired interface position on the boundary, such that $\cC\time{} \in \mathbb{R}^{2 \times \nstate}$, while 
\controlx[3]\time{}, in the third setting, uses seven outputs, i.e.\ $\cC\time{} \in \mathbb{R}^{7 \times \nstate}$.
With this, we compute the feedback gain matrices according to \cref{sec:LQR} with \cref{algo:BDF} and set $\ord = 1$.
Then, we simulate the closed-loop system with these feedback gain matrices together with the perturbations. 
The resulting feedback controls \controlx[1]\time{}, \controlx[2]\time{}, 
and \controlx[3]\time{} are displayed in the top part of \cref{fig:weights} and
the interface positions in the middle part of \cref{fig:weights}.
The interface positions are relative to the desired interface position
at the point $\statexstar\time$ with the largest deviation on the interface:
\begin{align*}
  \Gammaintdiff &= \GammaintRef - \Gammaint(\ttime, \statexstar), \\
  \statexstar &= \argmax_{\xx \in [0,0.5]}\ \left|\GammaintRefx - \Gammaint(\ttime, \xx)\right|.
\end{align*}

The feedback control \controlx[1]\time{} is most active shortly after the
perturbation starts and it steers the interface back to the desired position in
much reduced time, as expected from theory. 
Again as expected, a smaller weight factor \weight{} allows the control
\controlx[2]\time{} to react with a larger magnitude input while the perturbation is active.
It is, therefore, able to stop the interface from deviating earlier and drives
it back even faster.

It should be mentioned that smaller weights reduce the contribution of the control cost term
in the cost functional \labelcref{eq:cost}.
This term can also be interpreted as a regularization term. 
Consequently, smaller weights decrease the regularity of the feedback control problem.
We note that this results in larger computational cost for the solution of the DRE due to slower convergence
of the internal iterative solvers used in the BDF method, but refer to~\cite[Section~4.2]{BarBSS21a}
for a more detailed discussion of this issue.

The weights for the LQR settings of \controlx[1]\time{} and \controlx[2]\time{}
are not directly comparable to the setting for \controlx[3]\time, since this is based on different outputs.
The outputs at \Gammaoutx[5]\time{} and \Gammaoutx[6]\time{} allow it to detect the temperature perturbation earlier
and the output that monitors the movement of \Gammaint\time{} can observe the deviation of the interface earlier.
Thus, \controlx[3]\time{} is most active immediately after the perturbation starts and moves the interface back much faster.

Note, that we use a time-adaptive fractional-step-theta scheme
with a relative feedback control-based indicator function
to simulate the Stefan problem together with a feedback control.
For this example, the bottom part of \cref{fig:weights}
shows the adaptive time-step sizes.
While the feedback control is inactive, like at the beginning or end of the simulation, 
the method chooses the largest time-step size of $\timestepk = 0.0025$. 
As soon as the feedback control is active, the method chooses the minimal time-step size of $\timestepk = 10^{-4}$.
This prevents numerical instabilities that can occur with very small \weight{} and thus very active feedback controls
as is demonstrated in~\cite[Section~4.2]{BarBSS21a}.

\nameref{sec:ex1} demonstrates the influence of the chosen weight parameter \weight{} in the cost functional (\cref{eq:cost}). 
Smaller weight parameters increase the impact of the deviation term in the cost functional indirectly by
decreasing the relevance of the control costs which are measured by the second term and scaled by \weight.
Further, the choice and number of outputs has a significant impact on the performance of the feedback control.

The next experiment demonstrates that our method works as well with increased time-dependence of the data
in the non-autonomous DRE.

\subsection{Experiment 2}\label{sec:ex2} 
In our second experiment, the desired interface trajectory is, again, a flat horizontal line moving from its
initial height at \(0.5\), this time, upward by \(0.1\).
This is a longer distance interface movement than in \nameref{sec:ex1} on the same time horizon and, in turn,
results in a stronger time-dependence of the matrices for the DRE.\@ 
The investigated LQR design is chosen as in the second case in \nameref{sec:ex1}. 
\begin{figure}[t]
  \includegraphics{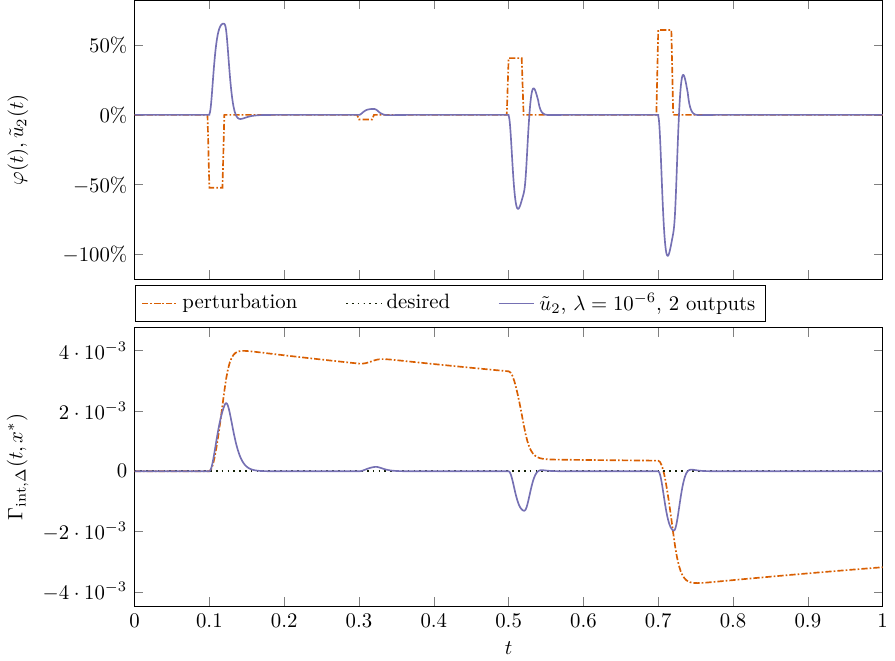}
  \caption{Perturbation and feedback (top) and relative interface position (bottom) for an interface moving upwards
	    (\nameref{sec:ex2})}%
  \label{fig:up}
\end{figure}

That means, we use the same weight ($\weight = 10^{-6}$) and outputs as for \controlx[2]\time{} 
above and call  the corresponding feedback control \controltilde[2]\time.
This time we use four randomly generated perturbations applied at the times $0.1, 0.3, 0.5$, and $0.7$, which, together with
\controltilde[2]\time, can be found in the top part of \cref{fig:up}.
Again, the perturbations cause the interface to deviate from the desired trajectory.
The feedback control \controltilde[2]\time{} behaves similarly to the previous experiment.
It stops the interface from deviating and drives it back to the desired position,
as expected from the theory.
Additionally, the time-adaptivity behaves analogously to \nameref{sec:ex1}.

The experiment showcases the performance of our method with strongly time-varying coefficients in the DRE.
We implemented more experiments with desired interface trajectories that move
considerably.
Our feedback control showed the same performance for all of them. This demonstrates
that our solver can also cope with the stronger time dependence without difficulties.

\subsection{Experiment 3}\label{sec:ex3} 
\begin{figure}[tp]
  \includegraphics{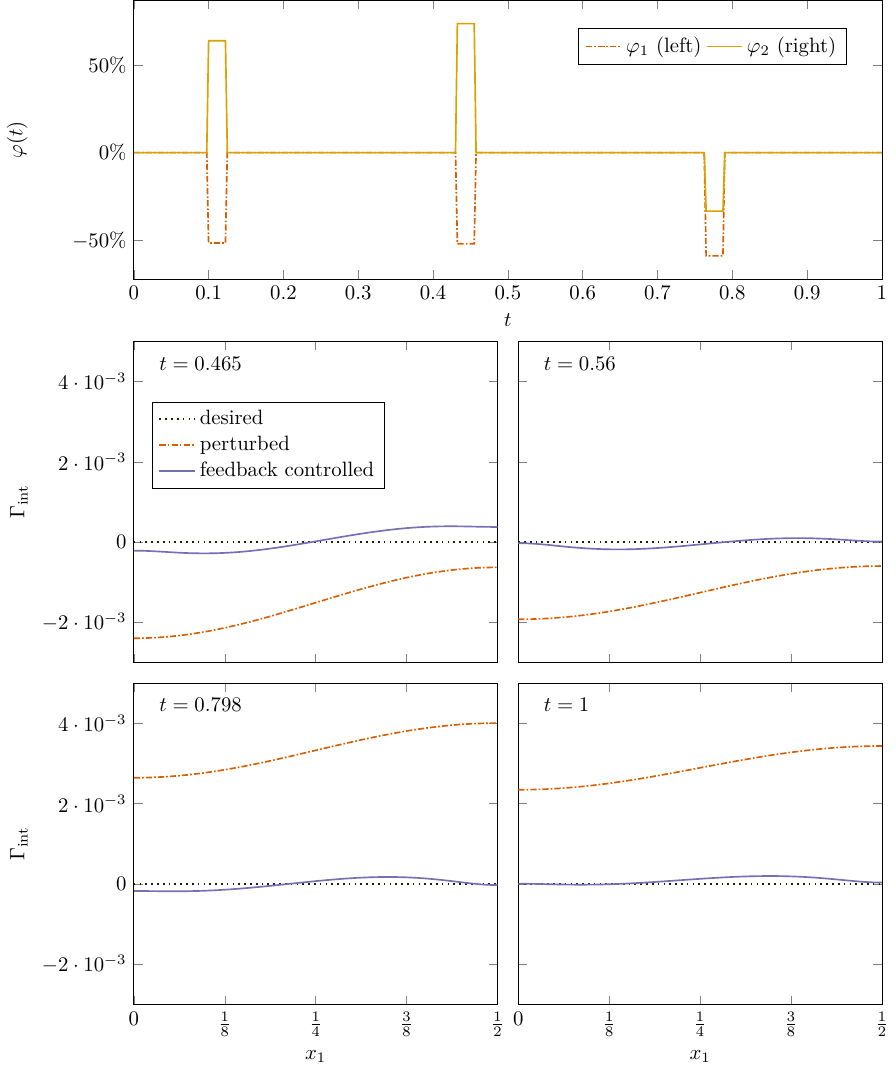}
  \caption{Two different perturbations (top), and perturbed and feedback controlled relative interface positions (bottom)
	    (\nameref{sec:ex3})}%
  \label{fig:curv}
\end{figure}

Our third experiment revisits the basic task as in \nameref{sec:ex1}. 
This time, we generate two different perturbations $\pertx[1]\time,
\pertx[2]\time \in [-\Tcool, \Tcool]$, acting on the left (\Gammacoolx[1]\time{}) and right (\Gammacoolx[2]\time{}) part
of the Dirichlet boundary at the bottom of \domain\time{} (see \cref{fig:inout}),
with three random values each, which are applied at the times $0.1, 0.4325$, and $0.765$. 
They are displayed in the top part of \cref{fig:curv}.
We apply these two perturbation functions to
\begin{align*}
  \temp	&= \Tcool + \pertx[1], \qquad \text{on}\ (0, \Endtime] \times \Gammacoolx[1],\\
  \temp	&= \Tcool + \pertx[2], \qquad \text{on}\ (0, \Endtime] \times \Gammacoolx[2].
\end{align*}
These perturbations will not only move the interface away from the desired position
but also add some curvature to the deviated interface.
To stabilize this interface position, we use six individual inputs at $\Gammainx[1]\time,\ldots,\Gammainx[6]\time$ 
($\cB\time{} \in \mathbb{R}^{\nstate \times 6}$) and the same two outputs at \Gammaoutx[3]\time{}, \Gammaoutx[4]\time{}
as in the previous experiments ($\cC\time{} \in \mathbb{R}^{2 \times \nstate}$).
For the weight factor in the cost functional we use $\weight = 10^{-9}$.
Again, this weight factor is not directly comparable to the previous experiments
since we use a different combination of inputs and outputs.

The resulting perturbed and feedback controlled interfaces
are displayed in the bottom part of \cref{fig:curv} at the
time points $\ttime \in [0.465, 0.56, 0.798, 1]$.
With the first two perturbations, the interface is pushed downwards and assumes a distinct curvature.
Right after the second perturbation ($\ttime = 0.465$), the feedback controlled interface is already essentially
back to the desired interface, but still has undesired curvature.
Well before the third perturbation, at $\ttime = 0.56$, the feedback controlled interface
is almost flat again, as the desired interface.
At the two endpoints of the interface, where the outputs are located, 
the controlled interface is directly at the desired position.
The remaining curvature can not be measured by the cost functional and is, thus, not 
actively removed.
For the given setup, this is the expected behavior of the feedback control.
It is similar after the first perturbation.
The third perturbation moves the interface upwards above the desired position ($\ttime = 0.798$).
At the end of the time interval ($\ttime = 1$), the feedback controlled interface is again almost flat and back
to the desired position at the ends while the uncontrolled interface is still far away.

This experiment shows that the feedback control can move the perturbed interface back to the desired trajectory
and additionally control the curvature of the interface, to a certain extent. 
With the chosen input and output setting, the position of the interface is corrected as fast as in the previous experiments.
The curvature correction requires some more time but is still performed during the considered time interval
up to the level that the cost functional can measure with the given two outputs.

\medskip
In all three experiments, we are able to stabilize the desired interface position
with the proposed feedback control approach
despite the different challenges that are posed by the Stefan problem.

\section{Conclusions}\label{sec:conclusions}
In this work, we address the challenging task to derive, compute, and apply a feedback control for the two-dimensional two-phase Stefan problem.
The particular challenges for the Stefan problem lie in the non-linearities, the DAE structure, as well as
the moving inner boundary.
Additionally, the applied feedback control approach results in a generalized DRE with time-dependent large-scale matrices.
We address this task with a sharp interface representation and mesh movement techniques.
Resulting from this coupling of the Stefan problem with mesh movement, we have a detailed description of the linearization, discretization and matrix assembly.
Regarding this, we particularly elaborate on how Dirichlet boundary conditions can be treated 
and display the DAE structure in the resulting matrices.

To obtain a feedback control for the Stefan problem, we apply the LQR approach and treat large-scale non-autonomous DREs with a corresponding low-rank BDF method.
More specifically, we also include the time-dependent mass matrix and its derivative.
The feedback control resulting from this approach is applied in a
forward simulation of the closed-loop Stefan problem
where we handle the numerical difficulties that arise with
a time-adaptive fractional-step-theta scheme specifically adapted to this process.

Through several numerical experiments, we demonstrate how effectively our methods and the resulting feedback controls perform.
The performance of the feedback controls strongly depends on the choice of control parameters like the weight factor in the cost functional
and the selected inputs and outputs.
The outputs are particularly important since in our formulation of the Stefan problem, the interface position to be controlled,
is not explicitly available as an output. 
Thus, the outputs need to indicate the interface deviation reliably.

Our methods include several stages of approximations, like the linearization, implicit index reduction, 
and low-rank representation of the DRE solution in an iterative solver.
Like this, we are able to make the numerical solution of this large-scale problem feasible
while the feedback controls still perform as known from problem types, that are already well studied in this regard.

Future research could investigate other promising methods for the non-autonomous DRE like, e.g., splitting schemes.
Additionally, a quantitative error analysis is still to be done in future work.

\bibliographystyle{alphaurl}
\bibliography{references}

\begin{thebibliography}{MOPP18}

\bibitem[AKFIJ03]{AboFIetal03}
H.~Abou-Kandil, G.~Freiling, V.~Ionescu, and G.~Jank.
\newblock {\em Matrix {R}iccati Equations in Control and Systems Theory}.
\newblock Systems \& Control: Foundations \& Applications. Birkh{\"a}user,
  Basel, Switzerland, 2003.
\newblock \href {https://doi.org/10.1007/978-3-0348-8081-7}
  {\path{doi:10.1007/978-3-0348-8081-7}}.

\bibitem[ANS14]{AntNS14}
H.~Antil, R.~H. Nochetto, and P.~Sodr\'e.
\newblock Optimal control of a free boundary problem: Analysis with
  second-order sufficient conditions.
\newblock {\em {SIAM} J. Control Optim.}, 52(5):2771--2799, 2014.
\newblock \href {https://doi.org/10.1137/120893306}
  {\path{doi:10.1137/120893306}}.

\bibitem[ANS15]{AntNS15}
H.~Antil, R.~H. Nochetto, and P.~Sodr\'e.
\newblock Optimal control of a free boundary problem with surface tension
  effects: A priori error analysis.
\newblock {\em {SIAM} J. Numer. Anal.}, 53(5):2279--2306, 2015.
\newblock \href {https://doi.org/10.1137/140958360}
  {\path{doi:10.1137/140958360}}.

\bibitem[AP98]{AscP98}
U.~M. Ascher and L.~R. Petzold.
\newblock {\em Computer Methods for Ordinary Differential Equations and
  Differential-Algebraic Equations}.
\newblock SIAM, Philadelphia, 1998.
\newblock \href {https://doi.org/10.1137/1.9781611971392}
  {\path{doi:10.1137/1.9781611971392}}.

\bibitem[Bar21]{LQR_Stefan_codes}
B.~Baran.
\newblock {Linear Quadratic Regulator Computation for a Stefan Problem using
  M.-M.E.S.S. and FEniCS}, April 2021.
\newblock \href {https://doi.org/10.5281/zenodo.4671886}
  {\path{doi:10.5281/zenodo.4671886}}.

\bibitem[BBH21]{BehBH21}
M.~Behr, P.~Benner, and J.~Heiland.
\newblock Invariant {G}alerkin trial spaces and {D}avison-{M}aki methods for
  the numerical solution of differential {R}iccati equations.
\newblock {\em Applied Mathematics and Computation}, 410:126401, 2021.
\newblock \href {https://doi.org/10.1016/j.amc.2021.126401}
  {\path{doi:10.1016/j.amc.2021.126401}}.

\bibitem[BBHS18]{BarBHetal18a}
B.~Baran, P.~Benner, J.~Heiland, and J.~Saak.
\newblock Optimal control of a {S}tefan problem fully coupled with
  incompressible {N}avier--{S}tokes equations and mesh movement.
\newblock {\em Analele Stiintifice ale Universitatii Ovidius Constanta: Seria
  Matematica}, XXVI(2):11--40, August 2018.
\newblock \href {https://doi.org/10.2478/auom-2018-0016}
  {\path{doi:10.2478/auom-2018-0016}}.

\bibitem[BBSS21]{BarBSS21a}
B.~Baran, P.~Benner, J.~Saak, and T.~Stillfjord.
\newblock Numerical methods for closed-loop systems with non-autonomous data.
\newblock e-print arXiv:X, arXiv, 2021.
\newblock math.NA.
\newblock URL: \url{https://arxiv.org/X}.

\bibitem[Ber10]{Ber10}
M.~Bernauer.
\newblock {\em Motion Planning for the Two-Phase {S}tefan Problem in Level Set
  Formulation}.
\newblock PhD thesis, Technische Universit\"at Chemnitz, Chemnitz, Germany,
  2010.
\newblock URL: \url{http://nbn-resolving.de/urn:nbn:de:bsz:ch1-qucosa-63654}.

\bibitem[BG91]{BitG91}
R.~R. Bitmead and M.~Gevers.
\newblock {\em Riccati Difference and Differential Equations: Convergence,
  Monotonicity and Stability}, pages 263--291.
\newblock Communications and Control Engineering. Springer Berlin Heidelberg,
  Berlin, Heidelberg, 1991.
\newblock \href {https://doi.org/10.1007/978-3-642-58223-3_10}
  {\path{doi:10.1007/978-3-642-58223-3_10}}.

\bibitem[BL18]{BenL18}
P.~Benner and N.~Lang.
\newblock Peer methods for the solution of large-scale differential matrix
  eqautions.
\newblock e-print arXiv:1804.08524, arXiv, July 2018.
\newblock math.NA.
\newblock URL: \url{https://arxiv.org/pdf/1804.08524.pdf}.

\bibitem[BLT09]{BenLT09}
P.~Benner, R.-C. Li, and N.~Truhar.
\newblock On the {ADI} method for {S}ylvester equations.
\newblock {\em J. Comput. Appl. Math.}, 233(4):1035--1045, 2009.
\newblock \href {https://doi.org/10.1016/j.cam.2009.08.108}
  {\path{doi:10.1016/j.cam.2009.08.108}}.

\bibitem[BM04]{BenM04}
P.~Benner and H.~Mena.
\newblock {BDF} methods for large-scale differential {R}iccati equations.
\newblock In B.~{De~Moor}, B.~Motmans, J.~Willems, P.~Van~Dooren, and
  V.~Blondel, editors, {\em Proc. 16th Intl. Symp. Mathematical Theory of
  Network and Systems, MTNS 2004}, 2004.

\bibitem[BM18]{BenM18}
P.~Benner and H.~Mena.
\newblock Numerical solution of the infinite-dimensional {LQR}-problem and the
  associated differential {R}iccati equations.
\newblock {\em J. Numer. Math.}, 26(1):1--20, March 2018.
\newblock published online May 2016.
\newblock \href {https://doi.org/10.1515/jnma-2016-1039}
  {\path{doi:10.1515/jnma-2016-1039}}.

\bibitem[BPS10]{BaePS10}
E.~B{\"a}nsch, J.~Paul, and A.~Schmidt.
\newblock An {ALE} {FEM} for solid-liquid phase transitions with free melt
  surface.
\newblock Technical report, Zentrum f{\"u}r Technomathematik, University of
  Bremen, 2010.
\newblock URL:
  \url{http://www.math.uni-bremen.de/zetem/cms/media.php/262/report1007.pdf}.

\bibitem[BPS13]{BaePS13}
E.~B{\"a}nsch, J.~Paul, and A.~Schmidt.
\newblock An {ALE} finite element method for a coupled {S}tefan problem and
  {N}avier–{S}tokes equations with free capillary surface.
\newblock {\em Internat. J. Numer. Methods Fluids}, 71(10):1282--1296, 2013.
\newblock \href {https://doi.org/10.1002/fld.3711}
  {\path{doi:10.1002/fld.3711}}.

\bibitem[FRM08]{morFreRM08}
F.~Freitas, J.~Rommes, and N.~Martins.
\newblock Gramian-based reduction method applied to large sparse power system
  descriptor models.
\newblock {\em {IEEE} Trans. Power Syst.}, 23(3):1258--1270, August 2008.
\newblock \href {https://doi.org/10.1109/TPWRS.2008.926693}
  {\path{doi:10.1109/TPWRS.2008.926693}}.

\bibitem[GHJK18]{GueHJetal18}
Y.~G{\"u}ldo{\v g}an, M.~Hached, K.~Jbilou, and M.~Kurulay.
\newblock Low rank approximate solutions to large-scale differential matrix
  {R}iccati equations.
\newblock {\em Applicationes Mathematicae}, 45(2):233--254, 2018.
\newblock \href {https://doi.org/10.4064/am2355-1-2018}
  {\path{doi:10.4064/am2355-1-2018}}.

\bibitem[Gup18]{Gup18}
S.~C. Gupta.
\newblock {\em The Classical Stefan Problem (Second Edition)}.
\newblock Elsevier, Amsterdam, 2018.
\newblock \href {https://doi.org/10.1016/B978-0-444-63581-5.09985-1}
  {\path{doi:10.1016/B978-0-444-63581-5.09985-1}}.

\bibitem[KK19]{KogK19}
S.~Koga and M.~Krstic.
\newblock Control of two-phase {S}tefan problem via single boundary heat input.
\newblock In {\em Proceedings of the IEEE Conference on Decision and Control},
  pages 2914--2919, 2019.
\newblock \href {https://doi.org/10.1109/CDC.2018.8619638}
  {\path{doi:10.1109/CDC.2018.8619638}}.

\bibitem[KK20a]{KogK20}
S.~Koga and M.~Krstic.
\newblock Single-boundary control of the two-phase {S}tefan system.
\newblock {\em Systems Control Lett.}, 135:104573, 2020.
\newblock \href {https://doi.org/10.1016/j.sysconle.2019.104573}
  {\path{doi:10.1016/j.sysconle.2019.104573}}.

\bibitem[KK20b]{KogK20a}
S.~Koga and M.~Krstic.
\newblock {\em Two-Phase Stefan Problem}, pages 139--157.
\newblock Springer International Publishing, Cham, 2020.
\newblock \href {https://doi.org/10.1007/978-3-030-58490-0_5}
  {\path{doi:10.1007/978-3-030-58490-0_5}}.

\bibitem[KM90a]{KunM90}
P.~Kunkel and V.~Mehrmann.
\newblock Numerical solution of differential algebraic {R}iccati equations.
\newblock {\em Linear Algebra Appl.}, 137/138:39--66, 1990.
\newblock \href {https://doi.org/10.1016/0024-3795(90)90126-W}
  {\path{doi:10.1016/0024-3795(90)90126-W}}.

\bibitem[KM90b]{KunM89}
P.~Kunkel and V.~Mehrmann.
\newblock Numerical solution of {R}iccati differential algebraic equations.
\newblock In M.~A.~Kaashoek et~al, editor, {\em Proceedings of the
  International Symposium on the Mathematical Theory of Networks and Systems,
  Amsterdam, Netherlands, June 1989}, pages 479--487, Basel, 1990.
  Birkh{\"a}user.

\bibitem[KM06]{KunM06}
P.~Kunkel and V.~Mehrmann.
\newblock {\em Differential-Algebraic Equations: Analysis and Numerical
  Solution}.
\newblock Textbooks in Mathematics. EMS Publishing House, Z{\"u}rich,
  Switzerland, 2006.

\bibitem[KM20]{KosM20}
A.~Koskela and H.~Mena.
\newblock Analysis of {K}rylov subspace approximation to large-scale
  differential {R}iccati equations.
\newblock {\em Electron. Trans. Numer. Anal.}, 52:431--454, 2020.
\newblock \href {https://doi.org/10.1553/etna_vol52s431}
  {\path{doi:10.1553/etna_vol52s431}}.

\bibitem[KS20]{KirSi19}
G.~Kirsten and V.~Simoncini.
\newblock Order reduction methods for solving large-scale differential matrix
  {R}iccati equations.
\newblock {\em {SIAM} J. Sci. Comput.}, 42(4):A2182--A2205, 2020.
\newblock \href {https://doi.org/10.1137/19M1264217}
  {\path{doi:10.1137/19M1264217}}.

\bibitem[Lan17]{Lan17}
N.~Lang.
\newblock {\em Numerical Methods for Large-Scale Linear Time-Varying Control
  Systems and related Differential Matrix Equations}.
\newblock {D}issertation, Technische Universit{\"a}t Chemnitz, Germany, June
  2017.
\newblock Logos-Verlag, Berlin, ISBN 978-3-8325-4700-4.
\newblock URL: \url{https://www.logos-verlag.de/cgi-bin/buch/isbn/4700}.

\bibitem[LC31]{LamC31}
G.~Lam{\'e} and B.~P. Clapeyron.
\newblock M{\'e}moire sur la solidification par refroidissement d’un globe
  liquide.
\newblock In {\em Annales Chimie Physique}, volume~47, pages 250--256, 1831.

\bibitem[LMS15]{LanMS15}
N.~Lang, H.~Mena, and J.~Saak.
\newblock On the benefits of the {$LDL^T$} factorization for large-scale
  differential matrix equation solvers.
\newblock {\em Linear Algebra Appl.}, 480:44--71, 2015.
\newblock \href {https://doi.org/10.1016/j.laa.2015.04.006}
  {\path{doi:10.1016/j.laa.2015.04.006}}.

\bibitem[LWH12]{dolfin}
A.~Logg, G.~N. Wells, and J.~Hake.
\newblock {\em DOLFIN: a C++/Python Finite Element Library}, chapter~10.
\newblock Springer, 2012.
\newblock \href {https://doi.org/10.1145/1731022.1731030}
  {\path{doi:10.1145/1731022.1731030}}.

\bibitem[LZL20]{LiZL20}
D.~Li, X.~Zhang, and R.~Liu.
\newblock Exponential integrators for large-scale stiff {R}iccati differential
  equations.
\newblock {\em J. Comput. Appl. Math.}, 389:113360, 2020.
\newblock \href {https://doi.org/10.1016/j.cam.2020.113360}
  {\path{doi:10.1016/j.cam.2020.113360}}.

\bibitem[Meh91]{Meh91}
V.~Mehrmann.
\newblock {\em The Autonomous Linear Quadratic Control Problem, Theory and
  Numerical Solution}.
\newblock Number 163 in Lecture Notes in Control and Information Sciences.
  Heidelberg, July 1991.

\bibitem[Men12]{Men12}
H.~Mena.
\newblock {\em Numerical Solution of Differential {R}iccati Equations Arising
  in Optimal Control Problems for Parabolic Partial Differential Equations}.
\newblock Unidad de Publicaciones de la Facultad de Ciencias, Quito-Ecuador,
  first edition, May 2012.
\newblock Available as ISBN: 978-9978-383-09-4.

\bibitem[MOPP18]{MenOPetal18}
H.~Mena, A.~Ostermann, L.-M. Pfurtscheller, and C.~Piazzola.
\newblock Numerical low-rank approximation of matrix differential equations.
\newblock {\em J. Comput. Appl. Math.}, 340:602--614, 2018.
\newblock \href {https://doi.org/10.1016/j.cam.2018.01.035}
  {\path{doi:10.1016/j.cam.2018.01.035}}.

\bibitem[NCM11]{NieCM11}
M.~Niezgodka, A.~Crowley, and A.~M. Meirmanov.
\newblock {\em The Stefan Problem}.
\newblock De Gruyter, 2011.
\newblock \href {https://doi.org/10.1515/9783110846720.245}
  {\path{doi:10.1515/9783110846720.245}}.

\bibitem[OPW18]{OstPW18}
A.~Ostermann, C.~Piazzola, and H.~Walach.
\newblock Convergence of a low-rank {L}ie--{T}rotter splitting for stiff matrix
  differential equations.
\newblock e-print arXiv:1803.10473, arXiv, March 2018.
\newblock math.NA.
\newblock URL: \url{https://arxiv.org/abs/1803.10473}.

\bibitem[Rei72]{Rei72}
W.~T. Reid.
\newblock {\em {R}iccati Differential Equations}, volume~86 of {\em Mathematics
  in Science and Engineering}.
\newblock Academic Press, New York, 1972.

\bibitem[Rub71]{Rub71}
L.~I. Ruben\v{s}te\u{\i}n.
\newblock {\em The {S}tefan problem}, volume~27 of {\em Translations of
  Mathematical Monographs}.
\newblock American Mathematical Society, Providence, R.I., 1971.
\newblock Translated from the Russian by A.~D. Solomon.
\newblock \href {https://doi.org/10.1090/mmono/027}
  {\path{doi:10.1090/mmono/027}}.

\bibitem[SKB21]{SaaKB21-mmess-2.1}
J.~Saak, M.~K{\"{o}}hler, and P.~Benner.
\newblock {M-M.E.S.S.}-2.1 -- {T}he {M}atrix {E}quations {S}parse {S}olvers
  library, April 2021.
\newblock see also: \url{https://www.mpi-magdeburg.mpg.de/projects/mess}.
\newblock \href {https://doi.org/10.5281/zenodo.4719688}
  {\path{doi:10.5281/zenodo.4719688}}.

\bibitem[Son98]{Son98}
E.~D. Sontag.
\newblock {\em Mathematical Control Theory}.
\newblock Texts in Applied Mathematics. Springer-Verlag, New York, NY, 2nd
  edition, 1998.
\newblock \href {https://doi.org/10.1007/978-1-4612-0577-7}
  {\path{doi:10.1007/978-1-4612-0577-7}}.

\bibitem[Ste89]{Ste89a}
J.~Stefan.
\newblock {\"U}ber einige {P}robleme der {T}heorie der {W}{\"a}rmeleitung.
\newblock {\em Sitzungber., Wien, Akad. Mat. Natur}, 98:473--484, 1889.

\bibitem[Ste90]{Ste90}
J.~Stefan.
\newblock {\"U}ber die {T}heorie der {E}isbildung.
\newblock {\em Monatshefte f{\"u}r Mathematik}, 1(1):1--6, 1890.

\bibitem[Ste91]{Ste89b}
J.~Stefan.
\newblock {\"U}ber die {T}heorie der {E}isbildung, insbesondere {\"u}ber die
  {E}isbildung im {P}olarmeere.
\newblock {\em Annalen der Physik und Chemie}, 42:269--286, 1891.

\bibitem[Sti15a]{Sti15}
T.~Stillfjord.
\newblock Low-rank second-order splitting of large-scale differential {R}iccati
  equations.
\newblock {\em {IEEE} Trans. Autom. Control}, 60(10):2791--2796, 2015.
\newblock \href {https://doi.org/10.1109/TAC.2015.2398889}
  {\path{doi:10.1109/TAC.2015.2398889}}.

\bibitem[Sti15b]{Sti15a}
T.~Stillfjord.
\newblock {\em Splitting schemes for nonlinear parabolic problems}.
\newblock PhD thesis, Lund University, 2015.
\newblock URL: \url{https://lup.lub.lu.se/search/ws/files/3905802/5277358.pdf}.

\bibitem[Sti18a]{Sti18a}
T.~Stillfjord.
\newblock Adaptive high-order splitting schemes for large-scale differential
  {R}iccati equations.
\newblock {\em Numer. Algorithms}, 78:1129--1151, 2018.
\newblock \href {https://doi.org/10.1007/s11075-017-0416-8}
  {\path{doi:10.1007/s11075-017-0416-8}}.

\bibitem[Sti18b]{Sti18}
T.~Stillfjord.
\newblock Singular value decay of operator-valued differential {L}yapunov and
  {R}iccati equations.
\newblock {\em {SIAM} J. Control Optim.}, 56:3598--3618, 2018.
\newblock \href {https://doi.org/10.1137/18M1178815}
  {\path{doi:10.1137/18M1178815}}.

\bibitem[Wei16]{Wei16}
H.~K. Weichelt.
\newblock {\em Numerical Aspects of Flow Stabilization by {R}iccati Feedback}.
\newblock {D}issertation, Otto-von-Guericke-Universit{\"a}t, Magdeburg,
  Germany, January 2016.
\newblock URL: \url{http://nbn-resolving.de/urn:nbn:de:gbv:ma9:1-8693}.

\bibitem[Zie08]{Zie08}
S.~Ziegenbalg.
\newblock {\em {K}ontrolle freier {R}{\"a}nder bei der {E}rstarrung von
  {K}ristallschmelzen}.
\newblock PhD thesis, {Technische Universit\"at Dresden}, Dresden, Germany,
  2008.
\newblock In German.
\newblock URL:
  \url{http://nbn-resolving.de/urn:nbn:de:bsz:14-ds-1212521184972-55836}.

\end{thebibliography}

\end{document}

